\documentclass[11pt, a4paper]{article}%
\usepackage{eurosym}
\usepackage{amsmath}
\usepackage{amsfonts}
\usepackage{bm}
\usepackage{amsfonts, graphicx, verbatim, amsmath,amssymb}
\usepackage{color}
\usepackage{amssymb}
\usepackage{graphicx}%
\providecommand{\U}[1]{\protect\rule{.1in}{.1in}}
\providecommand{\U}[1]{\protect\rule{.1in}{.1in}}
\newtheorem{theorem}{Theorem}[section]

\newtheorem{corollary}[theorem]{Corollary}
\newtheorem{corol}[theorem]{Corollary}
\newtheorem{Fact}[theorem]{Fact}

\newtheorem{lemma}[theorem]{Lemma}

\newtheorem{fact}[theorem]{Fact}
\newtheorem{lma}[theorem]{Lemma}

\newtheorem{proposition}[theorem]{Proposition}

\newtheorem{remark}[theorem]{Remark}
\newtheorem{Comment}[theorem]{Comment}

\newcommand{\rE}{\mathbb{E}}
\newcommand{\A}{{\mathcal{A}}}
\newcommand{\rP}{\mathbb{P}}
\newcommand{\p}{{\mathbb{P}}}
\newcommand{\Z}{{\mathbb{Z}}}

\newcommand{\E}{{\mathbb{E}}}

\newcommand{\eqdef}{\stackrel{\mathrm{def}}{=}}

\newcommand{\bb}{\begin{eqnarray*}}
\newcommand{\ee}{\end{eqnarray*}}
\newcommand{\bbb}{\begin{eqnarray}}
\newcommand{\eee}{\end{eqnarray}}

\newcommand{\qed}{$\diamond$}
\begin{document}

\title{Functional CLT for martingale-like nonstationary dependent structures }
\author{Florence Merlev\`ede{\thanks{Universit\'{e} Paris-Est, LAMA (UMR 8050), UPEM,
CNRS, UPEC. Email: florence.merlevede@u-pem.fr}}, Magda
Peligrad{\thanks{University of Cincinnati. Email: peligrm@ucmail.uc.edu}} and
Sergey Utev{\thanks{University of Leicester. Email: su35@leicester.ac.uk}}}
\maketitle

{\abstract} In this paper we develop non-stationary martingale techniques for
dependent data. We shall stress the non-stationary version of the projective
Maxwell-Woodroofe condition, which will be essential for obtaining maximal
inequalities and functional central limit theorem for the following examples:
nonstationary $\rho$-mixing sequences, functions of linear processes with
non-stationary innovations, quenched version of the functional central limit
theorem for a stationary sequence,  evolutions in random media such as a
process sampled by a shifted Markov chain.

\smallskip


\smallskip


\section{Introduction}

Historically, one of the most celebrated limit theorems in non-stationary
setting is, among others, the functional central limit theorem for
non-stationary sequences of martingale differences. For more general dependent
sequences, one of the basic techniques is approximate them with martingales by
using projection operators. A remarkable early result obtained by using this
technique is due to Dobrushin (1956), who studied non-stationary Markov
chains. Later the technique was used, also for Markov chains, in Sethuraman
and Varadhan (2005) and in Peligrad (2012). In order to treat more general
dependent structures, McLeish (1975, 1977) introduced the notion of
mixingales, which are martingale-like structures involving conditions imposed
to the bounds of the moments of projections of an individual variable on past
sigma fields. This method is very fruitful, but still involves a large degree
of stationarity and complicated additional assumptions. In general, the theory
of non-stationary martingale approximation it is much more difficult and it
has remained much behind the theory of martingale methods for stationary
processes. In the stationary setting, the theory of martingale approximations
was steadily developed. We mention the well-known results, such as the
celebrated results by Gordin (1969), Heyde (1974), Maxwell and Woodroofe
(2000) and newer results by Peligrad and Utev (2005), Zhao and Woodroofe
(2008), Gordin and Peligrad (2011), among many others. Inspired by these ideas
and using a direct martingale approach, we derive alternative conditions to
the mixingale-type conditions imposed by McLeish. Our projective conditions
lead to a non-stationary version of the weak invariance principle under the
so-called Maxwell-Woodroofe condition, which is known to be very sharp.
Surprisingly, also, is the fact that our approach leads directly to the
quenched invariance principle under the Maxwell-Woodroofe condition which was
first obtained by Cuny and Merlev\`{e}de (2014) with a completely different
proof. In addition, our approach is also efficient enough to get the
functional version of the central limit theorem for $\rho$-mixing sequences
satisfying the Lindeberg condition established in Utev (1990). For this class,
we completely answer an open problem raised by Ibragimov in 1991. Other
applications we shall consider are functions of linear processes with
nonstationary innovations and evolutions in random media, such as a process
sampled by a shifted Markov chain.

We begin by treating nonstationary sequences with the near linear behaviour of
the variance of the partial sums. Then, we disscuss the general non-stationary
triangular arrays and give the applications. The proofs are given in Section \ref{sectionproofs}.

\section{Results under the normalization $\sqrt{n}$}

Let $(X_{k})_{k\geq1}$ be a sequence of centered real-valued random variables
in ${\mathbb{L}}^{2}(\Omega,{\mathcal{A}},{\mathbb{P}})$ and set $S_{n}%
=\sum_{i=1}^{n}X_{i}$ for $n\geq1$ and $S_{0}=0$. Let $({\mathcal{F}}%
_{i})_{i\geq0}$ be a non-decreasing sequence of $\sigma$-algebras such that
$X_{i}$ is ${\mathcal{F}}_{i}$-measurable for any $i\geq1$. The following
notation will be often used: ${\mathbb{E}}_{k}(X):={\mathbb{E}}(X|{\mathcal{F}%
}_{k})$. In the sequel we denote by $D([0,1])$ the space of functions defined
on $[0,1]$, right continuous, with finite left hand limits, which is endowed
with uniform topology and by $[x]$ the integer part of $x$. For any $k\geq0$
let
\begin{equation}
\delta(k)=\max_{i\geq0}\Vert{\mathbb{E}}(S_{k+i}-S_{i}|{\mathcal{F}}_{i}%
)\Vert_{2}\,,\label{defakNS}%
\end{equation}
and for any $k,m\geq0$ let
\[
\theta_{k}^{m}=m^{-1}\sum_{i=1}^{m-1}{\mathbb{E}}_{k}(S_{k+i}-S_{k}) \, .
\]
To get the functional form of the central limit theorem under the
normalization $\sqrt{n}$, we shall assume the Lindeberg-type condition in the
form
\begin{equation}
\sup_{n\geq1}\frac{1}{n}\sum_{j=1}^{n}{\mathbb{E}}(X_{j}^{2})\leq
C<\infty,\mbox{ and }\lim_{n\rightarrow\infty}\frac{1}{n}\sum_{k=1}%
^{n}{\mathbb{E}}\{X_{k}^{2}I(|X_{k}|>\varepsilon\sqrt{n})\}=0\,,\,\text{ for
any $\varepsilon>0$.}\label{condLNS}%
\end{equation}

Our first result is in the spirit of Theorem 2.4 in McLeish (1977) and gives
sufficient conditions to ensure that the partial sums behave asymptotically
like a martingale.

\smallskip

\begin{theorem}
\label{CorollaryCLTNSforappli} Assume that the Lindeberg-type condition
(\ref{condLNS}) holds. Suppose in addition that
\begin{equation}
\sum_{k\geq0}2^{-k/2}\delta(2^{k})<\infty\, \label{MWcondNScasedyadic}%
\end{equation}
and there exists a constant $c^{2}$ such that, for any $t\in\lbrack0,1]$ and
any $\varepsilon>0$,
\begin{equation}
\lim_{m\rightarrow\infty}\limsup_{n\rightarrow\infty}{\mathbb{P}%
}\Big (\Big |\frac{1}{n}\sum_{k=1}^{[nt]}\big (X_{k}^{2}+2X_{k}\theta_{k}%
^{m}\big )-tc^{2}\Big |>\varepsilon\Big )=0\,. \label{condSqNS}%
\end{equation}
Then $\big \{n^{-1/2}\sum_{k=1}^{[nt]}X_{k},t\in\lbrack0,1]\big \}$ converges
in distribution in $D([0,1])$ to $cW$ where $W$ is a standard Brownian motion.
\end{theorem}

\begin{remark}
\label{remarkreinforced} Note that by the subadditivity property of the
sequence $(\delta(k))_{k\geq0}$, condition \eqref{MWcondNScasedyadic} is
equivalent to
\begin{equation}
\sum_{k\geq1}k^{-3/2}\delta(k)<\infty\,. \label{MWcondNScasenondyadic}%
\end{equation}
Moreover condition \eqref{MWcondNScasedyadic} holds under the stronger
condition
\begin{equation}
\sum_{k\geq1}k^{-1/2}\sup_{i\geq0}\Vert{\mathbb{E}}(X_{k+i}|{\mathcal{F}}%
_{i})\Vert_{2}<\infty\,. \label{MWcondNScasereinforced}%
\end{equation}

\end{remark}

For the stationary sequences, as a corollary to Theorem
\ref{CorollaryCLTNSforappli} we obtain:

\begin{corollary}
\label{rem2} Let $(X_{n})_{n\in{\mathbb{Z}}}$ be an ergodic stationary
sequence of centered random variables with finite second moment, which is
adapted to a stationary filtration $({\mathcal{F}}_{n})_{n\in{\mathbb{Z}}}$.
Assume that
\begin{equation}
\sum_{k\geq0}2^{-k/2}\Vert{\mathbb{E}}_{0}(S_{2^{k}})\Vert_{2}<\infty
\,,\label{MWcond}%
\end{equation}
Then, $\lim_{m\rightarrow\infty}m^{-1}\mathbb{E}(S_{m}^{2})=c^{2}$ and the
conclusion of Theorem \ref{CorollaryCLTNSforappli} holds.
\end{corollary}

Note that condition (\ref{MWcond}) is equivalent to $\sum_{k\geq1}%
k^{-3/2}\Vert{\mathbb{E}}_{0}(S_{k})\Vert_{2}<\infty$ and known under the name
of Maxwell-Woodroofe condition. Under this condition, Maxwell and Woodroofe
\cite{MW} obtained a CLT. Later, Peligrad and Utev \cite{PU1} have shown that
this condition is, in some sense, minimal in order for the sequence
$(S_{n}/\sqrt{n})_{n\geq1}$ to be stochastically bounded and they proved a
maximal inequality and convergence to the Brownian motion. In order to derive
this corollary from Theorem \ref{CorollaryCLTNSforappli} we use the fact that
$\delta(k)=\Vert{\mathbb{E}}(S_{k}|{\mathcal{F}}_{0})\Vert_{2}$ and then
condition (\ref{MWcondNScasedyadic}) reads as condition (\ref{MWcond}). In
addition, for $k\geq0$, we get, by the ergodic theorem,
\[
\lim_{n\rightarrow\infty}{\mathbb{E}}\Big \vert\frac{1}{n}\sum_{k=1}%
^{[nt]}(X_{k}^{2}+2X_{k}\theta_{k}^{m})-c^{2}%
t\Big \vert=t\big \vert{\mathbb{E}}X_{0}^{2}+2{\mathbb{E}}(X_{0}\theta_{0}%
^{m}\big )-c^{2}\big \vert\,.
\]
It remains to take into account that
\[
\frac{1}{m}\mathbb{E}(S_{m}^{2})={\mathbb{E}}(X_{0}^{2})+\frac{2}{m}\sum
_{i=1}^{m-1}\sum_{j=1}^{m-i}{\mathbb{E}}(X_{0}X_{j})={\mathbb{E}}(X_{0}%
^{2})+2{\mathbb{E}}(X_{0}\theta_{0}^{m}\big )\,,
\]
proving the corollary since it has been shown in \cite{PU1} that, in the
stationary setting, condition \eqref{MWcond} implies that $\lim_{m\rightarrow
\infty}m^{-1}\mathbb{E}(S_{m}^{2})$ exists.

\section{Results for general triangular arrays}

Let $\{X_{i,n},1\leq i\leq n\}$ be a triangular array of of square integrable
(${\mathbb{E}}(X_{i,n}^{2})<\infty$), centered (${\mathbb{E}}(X_{i,n})=0$),
real-valued random variables adapted to a filtration $(\mathcal{F}%
_{i,n})_{i\geq0}$. We write as before ${\mathbb{E}}_{j,n}(X)={\mathbb{E}%
}(X|\mathcal{F}_{j,n})$ and set 
\[
S_{k,n}=\sum_{i=1}^{k}X_{i,n} \, \text{ and } \, \theta_{k,n}^{m}=m^{-1}\sum_{i=1}^{m-1}{\mathbb{E}}_{k,n}(S_{k+i,n}-S_{k,n}) \, .
\]
We shall
assume that the triangular array satisfies the following Lindeberg type
condition:
\begin{equation}
\sup_{n\geq1}\sum_{j=1}^{n}{\mathbb{E}}(X_{j,n}^{2})\leq C<\infty,\mbox{
and }\lim_{n\rightarrow\infty}\sum_{k=1}^{n}{\mathbb{E}}\{X_{k,n}%
^{2}I(|X_{k,n}|>\varepsilon)\}=0\,,\,\text{ for any $\varepsilon>0$%
.}\label{condLNS_g}%
\end{equation}

Moreover, for a non-negative integer $u$ and positive integers $\ell,m$,
define martingale-type dependence characteristics by
\begin{equation}
\label{defalu_g}A^{2}(u)=\sup_{n\geq1}\sum_{k=0}^{n-1}\Vert{\mathbb{E}}%
_{k,n}(S_{k+u,n}-S_{k,n})\Vert_{2}^{2}\,.
\end{equation}
and
\begin{equation}
\label{defblu_g}B^{2}(\ell,m)=\sup_{n\geq1}\sum_{k=0}^{[n/\ell]}\Vert{\bar{S}%
}_{k,n}(\ell,m)\Vert_{2}^{2}\,,
\end{equation}
where
\[
{\bar{S}}_{k,n}(\ell,m)=\frac{1}{m}\sum_{u=0}^{m-1}\big ({\mathbb{E}%
}_{(k-1)\ell+1,n}(S_{(k+1){\ell}+u,n}-S_{k\ell+u,n})\big )\,.
\]

Our next theorem provides a general functional CLT under the Lindeberg
condition for martingale-like nonstationary triangular arrays.

\begin{theorem}
\label{fclt_MPU} Assume that the Lindeberg-type condition \eqref{condLNS_g}
holds and that
\begin{equation}
\lim_{j\rightarrow\infty}2^{-j/2}A(2^{j})=0\text{ and }\liminf_{j\rightarrow
\infty}\sum_{\ \ell\geq j}B(2^{\ell},2^{j})=0\,. \label{MWcondNS2_g}%
\end{equation}
Moreover, assume in addition that there exist a sequence of non-decreasing and
right-continuous functions $v_{n}(\cdot):[0,1]\rightarrow\{0,1,2,\ldots,n\}$
and a non-negative real $c^{2}$ such that, for any $t\in(0,1]$,
\begin{equation}
\lim_{m\rightarrow\infty}\limsup_{n\rightarrow\infty}{\mathbb{P}%
}\Big (\Big |\sum_{k=1}^{v_{n}(t)}\big (X_{k,n}^{2}+2X_{k,n}\theta_{k,n}%
^{m}\big )-tc^{2}\Big |>\varepsilon\Big )=0\, . \label{AssGr3_1_g}%
\end{equation} 
Then $\big \{\sum_{k=1}^{v_{n}(t)}X_{k,n},t\in\lbrack0,1]\big \}$ converges in
distribution in $D([0,1])$ to $cW$ where $W$ is a standard Brownian motion.
\end{theorem}

The following proposition is useful for verifying condition \eqref{AssGr3_1_g}.

\begin{proposition}
\label{toverifycondSqNS_g} Assume that the Lindeberg-type condition
\eqref{condLNS_g} holds. Assume in addition that for any non-negative integer
$\ell$,
\begin{equation}
\lim_{b\rightarrow\infty}\limsup_{n\rightarrow\infty}\sum_{k=b+1}^{n}%
\Vert{\mathbb{E}}_{k-b,n}(X_{k,n}X_{k+\ell,n})-{\mathbb{E}}_{0,n}%
(X_{k,n}X_{k+\ell,n})\Vert_{1}=0\,\label{convL1carreNSrho}%
\end{equation}
and, for any $t\in\lbrack0,1]$,
\begin{equation}
\lim_{m\rightarrow\infty}\limsup_{n\rightarrow\infty}{\mathbb{P}%
}\Big (\Big |\sum_{k=1}^{v_{n}(t)}\big ({\mathbb{E}}_{0,n}(X_{k,n}%
^{2})+2{\mathbb{E}}_{0,n}(X_{k,n}\theta_{k,n}^{m})\big )-tc^{2}%
\Big |>\varepsilon\Big )=0\,.\label{condSqNSavecE0rho}%
\end{equation}
Then the convergence \eqref{AssGr3_1_g} is satisfied.
\end{proposition}

\medskip

Let us apply the general Theorem \ref{fclt_MPU} to the sequences of random
variables when the normalizing sequence is $\sqrt{n}$. For any non-negative
integer $u$ and any positive integers $\ell$ and $m$, let $a(u)$ and
$b(\ell,m)$ be the non-negative quantities defined by
\begin{align*}
a^{2}(u) &  =\sup_{n\geq1}\frac{1}{n}\sum_{k=0}^{n-1}\Vert{\mathbb{E}}%
_{k}(S_{k+u}-S_{k})\Vert_{2}^{2}\quad,\quad b^{2}(\ell,m)=\sup_{n\geq1}%
\frac{1}{n}\sum_{k=0}^{n-1}\Vert{\bar{S}}_{k}(\ell,m)\Vert_{2}^{2}\,,\\
&  \mbox{where}\quad{\bar{S}}_{k}(\ell,m)=\frac{1}{m}\sum_{u=0}^{m-1}%
\big ({\mathbb{E}}_{(k-1)\ell+1}(S_{(k+1){\ell}+u}-S_{k\ell+u})\big )\,.
\end{align*}
The conditions are in the sense an average version of condition
\eqref{MWcondNScasedyadic}. They are particularly useful in the analysis of
quenched limit theorems. By applying Theorem \ref{fclt_MPU} to the triangular
array $X_{k,n}=X_{k}/\sqrt{n};$ $1\leq k\leq n$ and $v_{n}(t)=[nt]$ we obtain
the following corollary:

\begin{corollary}
\label{comment1} The statement of Theorem \ref{CorollaryCLTNSforappli} holds
when condition \eqref{MWcondNScasedyadic} is replaced by the following
conditions
\begin{equation}
\lim_{j\rightarrow\infty}2^{-j/2}a(2^{j})=0\text{ and }\liminf_{j\rightarrow
\infty}\sum_{\ell\geq j}2^{-\ell/2}b(2^{\ell},2^{j})=0\,. \label{MWcondNS2}%
\end{equation}

\end{corollary}

By using the definition of $(\delta(k))_{k\geq1}$ in (\ref{defakNS}),\ the
subadditivity of this sequence and Proposition 2.5 in Peligrad and Utev
(2015)\ we note that condition (\ref{MWcondNScasedyadic}) implies
$\lim_{j\rightarrow\infty}2^{-j/2}a(2^{j})=0.$ Moreover, condition
(\ref{MWcondNScasedyadic}) easily implies the second part of condition
(\ref{MWcondNS2}). By using this remark we can see that Theorem
\ref{CorollaryCLTNSforappli} is a consequence of Corollary \ref{comment1}. We
elected to present the results first for sequences of random variables and
then for triangular arrays, for stressing the fact that our results are
generalization to nonstationary sequences of the important results in the
stationary setting involving condition (\ref{MWcond}). The results are also
related to conditions in McLeish (1975, 1977). Our approach uses a suitable
martingale approximation whereas McLeish (1975, 1977) proved first tightness
of the partial sum process and then he identified the limit by using a
suitable characterization of the Wiener process given in Theorem 19.4 in
Billingsley (1968).

\section{Applications}

\subsection{$\rho$-mixing triangular arrays and sequences}

For a triangular array $\{X_{i,n},1\leq i\leq n\}$ of square integrable
(${\mathbb{E}}(X_{i,n}^{2})<\infty$), centered (${\mathbb{E}}(X_{i,n})=0$),
real-valued random variables, we denote by $\sigma_{k,n}^{2}=\mathrm{Var}%
\big (\sum_{\ell=1}^{k}X_{\ell,n}\big )$ for $k\leq n$ and $\sigma_{n}%
^{2}=\sigma_{n,n}^{2}$. For $0\leq t\leq1$, let
\[
v_{n}(t)=\inf\Big \{k;\,1\leq k\leq n\colon\frac{\sigma_{k,n}^{2}}{\sigma
_{n}^{2}}\geq t\Big \}\,\text{ and }\,W_{n}(t)=\sigma_{n}^{-1}\sum
_{i=1}^{v_{n}(t)}X_{i,n}\,.
\]
Define also $S_{k}=S_{k,n}=\sum_{i=1}^{k}X_{i,n}$. In this section we assume
that the triangular array is $\rho$-mixing in the sense that
\[
\rho(k)=\sup_{n\geq1}\max_{1\leq j\leq n-k}\rho\big (\sigma(X_{i,n},1\leq
i\leq j),\sigma(X_{i,n},j+k\leq i\leq n)\big )\rightarrow0\,,\,\text{as
$k\rightarrow\infty$}\,,
\]
where $\sigma(X_{t},t\in A)$ is the $\sigma$-field generated by the r.v.'s
$X_{t}$ with indices in $A$ and we recall that the maximal correlation
coefficient $\rho({\mathcal{U}},{\mathcal{V}})$ between two $\sigma$-algebras
is defined by
\[
\rho(\mathcal{U},\mathcal{V})=\sup\{|\mathrm{corr}(X,Y)|:X\in{\mathbb{L}}%
^{2}(\mathcal{U}),Y\in{\mathbb{L}}^{2}(\mathcal{V})\}\,.
\]
Next result gives the functional version of the central limit theorem for
$\rho$-mixing sequences satisfying the Lindeberg condition established in
Theorem 4.1 in Utev (1990). It answers an open question raised by Ibragimov in 1991.

\begin{theorem}
\label{ThrhoNSalter} Suppose that
\begin{equation}
\sup_{n\geq1}\sigma_{n}^{-2}\sum_{j=1}^{n}{\mathbb{E}}(X_{j,n}^{2})\leq
C<\infty\,,\label{condboundvar}%
\end{equation}
and
\begin{equation}
\lim_{n\rightarrow\infty}\sigma_{n}^{-2}\sum_{k=1}^{n}{\mathbb{E}}%
\{X_{k,n}^{2}I(|X_{k,n}|>\varepsilon\sigma_{n})\}=0\,,\,\text{ for any
$\varepsilon>0$}\,.\label{condLNSavecsigma}%
\end{equation}
Assume in addition that
\begin{equation}
\sum_{k\geq0}\rho(2^{k})<\infty\,.\label{condonrho}%
\end{equation}
Then $\big \{W_{n}(t),t\in(0,1]\big \}$ converges in distribution in
$D([0,1])$ (equipped with the uniform topology) to $W$ where $W$ is a standard
Brownian motion.
\end{theorem}

For the $\rho-$mixing sequences we also obtain the following corollary:

\begin{corollary}
\label{ThrhoNS} Let $(X_{n})_{n\geq1}$ be a sequence of centered random
variables in ${\mathbb{L}}^{2}({\mathbb{P}})$. Let $S_{n}=\sum_{k=1}^{n}X_{k}$
and $\sigma_{n}^{2}=\mathrm{Var}(S_{n})$. Suppose that conditions
\eqref{condboundvar}, \eqref{condLNSavecsigma} and \eqref{condonrho} are
satisfied. In addition assume that $\sigma_{n}^{2}=nh(n)$ where $h$ is a
slowly varying function at infinity. Then $\big \{\sigma_{n}^{-1}\sum
_{k=1}^{[nt]}X_{k},t\in(0,1]\big \}$ converges in distribution in $D([0,1])$
to $W$ where $W$ is a standard Brownian motion.
\end{corollary}

It should be noted that if $\big \{ \sigma_{n}^{-1} \sum_{k=1}^{[nt]} X_{k} ,
t \in[0,1] \big \}$ converges weakly to a standard Brownian motion, then
necessarily $\sigma_{n}^{2}=nh(n)$ where $h(n)$ is a slowly varying function
(i.e. a regularly varying function with exponent $1$). This is so since for $t
\in[0,1]$ fixed we have $S_{[nt]}/\sigma_{n}\rightarrow^{d} N(0,t)$ and in
addition, taking $t=1$ we have $S_{n}^{2}/\sigma_{n}^{2}$ is uniformly
integrable (by the convergence of moments theorem), implying $\sigma_{\lbrack
nt]}^{2}/\sigma_{n}^{2}\rightarrow t$.

\begin{Comment}
\textrm{If in Corollary \ref{ThrhoNS} above we assume that $\sigma_{n}^{2} =
n^{\alpha} h(n)$ where $\alpha>0$ and $h$ is a slowly varying function at
infinity, then the proof reveals that, under \eqref{condboundvar},
\eqref{condLNSavecsigma} and \eqref{condonrho}, $\big \{ \sigma_{n}^{-1}
\sum_{k=1}^{[nt]} X_{k} , t \in[0,1] \big \}$ converges in distribution in
$D([0,1])$ to $\{G(t), t \in[0,1]\} $ where, for any $t \in[0,1]$, $G(t) =
\sqrt{\alpha} \int_{0}^{t} u^{(\alpha-1) /2 } dW (u)$ with $W$ a standard
Brownian motion. }
\end{Comment}

\subsection{Functions of linear processes}

Let $(\varepsilon_{i})_{i\in{{\mathbb{Z}}}}$ be a sequence of real-valued
independent random variables. We shall say that the sequence $(\varepsilon
_{i})_{i\in{{\mathbb{Z}}}}$ satisfies the condition $(A)$ if $\,(\varepsilon
_{i}^{2})_{i\in{{\mathbb{Z}}}}$ is an uniformly integrable family and
$\sup_{i\in{\mathbb{Z}}}\Vert\varepsilon_{i}\Vert_{2}:=\sigma_{\varepsilon
}<\infty$. Let $(a_{i})_{i\geq0}$ be a sequence of reals in $\ell^{1}$. For
any integer $k$, let then $Y_{k}=\sum_{i\geq0}a_{i}\varepsilon_{k-i}$. Let $f$
be a function from ${\mathbb{R}}$ to ${\mathbb{R}}$ in the class
${\mathcal{L}}(c)$, meaning that there exists a concave non-decreasing
function $c$ from ${\mathbb{R}}^{+}$ to ${\mathbb{R}}^{+}$ with $\lim
_{x\rightarrow0}c(x)=0$ and such that
\[
|f(x)-f(y)|\leq c(|x-y|)\,\text{ for any }\,(x,y)\in{\mathbb{R}}^{2}\,.
\]
We shall also assume that
\begin{equation}
c\Big (K\sum_{i\geq0}|a_{i}|\Big )<\infty\,\text{ for any finite real $K>0$
and }\sum_{k\geq1}k^{-1/2}c\Big (2\sigma_{\varepsilon}\sum_{i\geq k}%
|a_{i}|\Big )<\infty\,,\label{Kac}%
\end{equation}
and, for any $k\geq1$, define
\begin{equation}
X_{k}=f(Y_{k})-{\mathbb{E}}(f(Y_{k}))\,.\label{deffunctionlinearXk}%
\end{equation}

\begin{corollary}
\label{corforlinear} Let $(\varepsilon_{i})_{i\in{{\mathbb{Z}}}}$ be a
sequence of real-valued independent random variables satisfying the condition
$(A)$. Let $f$ be a function from ${\mathbb{R}}$ to ${\mathbb{R}}$ belonging
to the class ${\mathcal{L}}(c)$ and let $(a_{i})_{i\geq0}$ be a sequence of
reals in $\ell^{1}$. Assume that condition \eqref{Kac} is satisfied and define
$(X_{k})_{k\geq1}$ by \eqref{deffunctionlinearXk}. Let $S_{n}=\sum_{k=1}%
^{n}X_{k}$ and $\sigma_{n}^{2}=\mathrm{Var}(S_{n})$. If $\sigma_{n}^{2}=nh(n)$
where $h(n)$ is a slowly varying function at infinity such that $\liminf
_{n\rightarrow\infty}h(n)>0$, then $\big \{\sigma_{n}^{-1}\sum_{k=1}%
^{[nt]}X_{k},t\in\lbrack0,1]\big \}$ converges in distribution in $D([0,1])$
to $W$ where $W$ is a standard Brownian motion.
\end{corollary}

Note that, if $|a_{i}|\leq C \rho^{i}$ for some $C>0$ and $\rho\in]0, 1[$, the
condition (\ref{Kac}) holds as soon as:
\[
\int_{0}^{1} \frac{c(t)}{t \sqrt{|\log t|}} dt < \infty\, .
\]
Note that this condition is satisfied as soon as $c(t)\leq D|\log
(t)|^{-\gamma}$ for some $D>0$ and some $\gamma> 1/2$. In particular, it is
satisfied if $f$ is $\alpha$-H\"older for some $\alpha\in]0, 1]$.

\begin{Comment}
\textrm{Note that other conditions on the innovations $(\varepsilon_{i})_{i
\in{{\mathbb{Z}}}}$ can lead to the conclusion of Corollary \ref{corforlinear}%
. For instance, let $(\phi_{\varepsilon}(k))_{k \geq0}$ be the sequence of
$\phi$-mixing coefficients associated with $(\varepsilon_{i})_{i
\in{{\mathbb{Z}}}}$. These coefficients are defined in the following way: for
any non-negative integer $k$,
\[
\phi_{\varepsilon}(k)= \sup_{\ell\in{\mathbb{Z}}} \phi({\mathcal{F}}_{\ell},
{\mathcal{G}}_{\ell+k } )
\]
with ${\mathcal{F}}_{\ell}= \sigma( \varepsilon_{i}, i \leq\ell)$,
${\mathcal{G}}_{\ell}= \sigma( \varepsilon_{i}, i \geq\ell)$ and where, for
two $\sigma$-algebras of ${\mathcal{A}}$,
\[
\phi(\mathcal{U},\mathcal{V}) =\sup\{ | {\mathbb{P}} ( V | U ) - {\mathbb{P}}
(V) | : U \in\mathcal{U}, {\mathbb{P}} (U) >0, V\in\mathcal{V} \} \, .
\]
Analyzing the proof of Corollary \ref{corforlinear} and taking into account
the computations made in Proposition 12 of \cite{DMPU}, we infer that
Corollary \ref{corforlinear} still holds if we impose the following conditions
on the innovations: $\sup_{i \in{\mathbb{Z}}} \Vert\varepsilon_{i}
\Vert_{\infty}:=B < \infty$ and $\sum_{k \geq0} \phi_{\varepsilon}(k) <
\infty$. }
\end{Comment}

\subsection{Quenched functional central limit theorems}

In this subsection we start with a stationary sequence and address the
question of functional CLT when the process is not started from its
equilibrium, but it is rather started at a point or from a fixed past
trajectory. This process is no longer strictly stationary. This type of result
is known under the name of quenched limit theorem. It is convenient to
introduce a stationary process by using the dynamical systems language. Let
$(\Omega,\mathcal{A}, {\mathbb{P}})$ be a probability space, and
$T:\Omega\mapsto\Omega$ be a bijective bimeasurable transformation preserving
the probability ${\mathbb{P}}$. An element $A$ is said to be invariant if
$T(A)=A$. We denote by $\mathcal{I}$ the $\sigma$-algebra of all invariant
sets. The probability ${\mathbb{P}}$ is ergodic if each element of
$\mathcal{I}$ has measure 0 or 1.

Let $\mathcal{F}_{0}$ be a $\sigma$-algebra of $\mathcal{A}$ satisfying
$\mathcal{F}_{0}\subseteq T^{-1}(\mathcal{F}_{0})$ and define the
nondecreasing filtration $(\mathcal{F}_{i})_{i\in{\mathbb{Z}}}$ by
$\mathcal{F}_{i}=T^{-i}(\mathcal{F}_{0})$. We assume that there exists a
regular version $P_{T|{\mathcal{F}}_{0}}$ of $T$ given ${\mathcal{F}}_{0}$, 

In this subsection, we assume that ${\mathbb{P}}$ is ergodic and we consider
$X_{0}$ a $\mathcal{F}_{0}$-measurable, square integrable and centered random
variable. Define then the sequence $\mathbf{X}=(X_{i})_{i\in\mathbb{Z}}$ by
$X_{i}=X_{0}\circ T^{i}$. Let $S_{n}=X_{1}+\cdots+X_{n}$ and $W_{n}%
=\{W_{n}(t),t\in\lbrack0,1]\}$ where $W_{n}(t)=n^{-1/2}S_{[nt]}$. It is
well-known that, by a canonical construction, any stationary sequence can be
represented in this way via the translation operator. As we shall see,
applying our Corollary \ref{comment1}, we derive the following quenched CLT in
its functional form under Maxwell and Woodroofe condition \eqref{MWcond}
which, from the subadditivity property of the sequence $(\Vert{\mathbb{E}}%
_{0}(S_{n})\Vert_{2})_{n\geq0}$, is equivalent to the convergence:
$\sum_{k\geq1}k^{-3/2}\Vert{\mathbb{E}}_{0}(S_{k})\Vert_{2}<\infty$. This
result was first obtained by Cuny and Merlev\`{e}de in 2014 (see their Theorem
2.7) with a completely different proof.

\begin{corollary}
\label{corquenched} Assume that \eqref{MWcond} holds. Then there exists a
constant $c^{2}$ such that $\lim_{n\rightarrow\infty}n^{-1/2}{\mathbb{E}%
}(S_{n}^{2})=c^{2}$ and $W_{n}$ satisfies the following quenched weak
invariance principle: on a set of probability one, for any continuous and
bounded function $f$ from $(D([0,1),\Vert\cdot\Vert_{\infty})$ to
${\mathbb{R}}$,
\[
\lim_{n\rightarrow\infty}{\mathbb{E}}_{0}(f(W_{n}))=\int f(zc)W(dz)\,,
\]
where $W$ is the distribution of a standard Wiener process.
\end{corollary}

The conclusion of this corollary can also be expressed in the following way.
Denote by ${\mathbb{P}^{\omega}}(A)$ a regular version of conditional
probability ${\mathbb{P}}(A|{\mathcal{F}}_{0})(\omega)$. Then for any $\omega$
in a set of probability $1$, $W_{n}$ converges in distribution in $D([0,1])$
to $W$ under ${\mathbb{P}^{\omega}}$.

Since condition (\ref{MWcond}) is verified by a stationary $\rho-$mixing
sequence satisfying (\ref{condonrho}) (see for instance Peligrad and Utev
\cite{PU1}) we can formulate the following Corollary.

\begin{corol}
\label{corquenchedro}Assume that condition (\ref{condonrho}) holds. Then the
quenched functional CLT in Corollary \ref{corquenched} holds.
\end{corol}

As a particular example we shall give now a quenched functional central limit
theorem for functions of a stationary Gaussian process.

Denote by $T$ the unit circle in the complex plane. Let $\mu$ denote
normalized Lebesgue measure on $T$ (normalized so that $\mu(T)=1$). For a
given stationary random sequence $X:=(X_{k})_{k\in{\mathbb{Z}}}$, a
\textquotedblleft spectral density function\textquotedblright\ (if one exists)
can also be viewed as a real, nonnegative, Borel, integrable function $f
:T\rightarrow\lbrack0,\infty)$ such that $f(\mathrm{e}^{i\lambda}) =
f(\mathrm{e}^{-i\lambda})$ for all $\lambda\in(0, \pi]$ and, for every
$k\in{\mathbb{Z}}$,
\[
\mathrm{cov}(X_{k},X_{0})=\int\limits_{t\in T}t^{k}f(t)\mu(dt) \, .
\]
Let $(a_{n})_{n\geq1}$ be a nonincreasing sequence of positive numbers such
that, $a_{n} \rightarrow0$, as $n \rightarrow\infty$, $a_{n}/a_{2n} \leq6/5$,
for all $n \geq1$, and $\sum_{k\geq0}a_{2^{k}}<\infty$. Let $Y:=(Y_{k}%
)_{k\in{\mathbb{Z}}}$ be a stationary Gaussian process with spectral density
of the form%
\[
f(\mathrm{e}^{i\lambda})=\exp\Big ( \sum\nolimits_{k=1}^{\infty}k^{-1}%
a_{k}\cos(k\lambda) \Big ) \, .
\]
For a fixed positive integer $M$ and a measurable function $g:{\mathbb{R}}%
^{M}\rightarrow{\mathbb{R}}$, defined the stationary sequence $(X_{k}%
)_{k\in{\mathbb{Z}}}$ by
\begin{equation}
X_{k}=g(Y_{k+1},Y_{2}, \ldots,Y_{k+M}) \, . \label{defX}%
\end{equation}

\begin{corol}
\label{Gaussian}Assume that $(X_{k})_{k\in{\mathbb{Z}}}$ defined by
(\ref{defX}) is a sequence of centered square integrable random variables.
Then the conclusion of Corollary \ref{corquenched} holds.
\end{corol}

\subsection{Application to a random walk in random time scenery}

Consider the partial sums associated with $(X_{k})_{k \geq0}$ which is  a
sequence of random variables, $\{\zeta_{j}\}_{j \geq0 }$, called the
\textit{random time scenery}, sampled by the process $(Y_{k})_{k\geq0}$,
defined as
\[
Y_{k} = k+ \phi_{k}, \quad k\geq0,
\]
where $\{\phi_{n}\}_{n\geq0}$ is a ``renewal"-type Markov chain defined as
follows: $\{\phi_{k};k\geq0\}$ is a discrete Markov chain with the state space
$\Z^+$ and transition matrix $P=(p_{ij})$ given by $p_{k,k-1}=1$
for $k\geq1$ and $p_{j}=p_{0,j-1}={\mathbb{P}}(\tau=j)$, $j=1,2,\ldots,$ (that
is whenever the chain hits $0$ it then regenerates with the
probability $p_{j}$). Therefore the sequence $(X_{k})_{k \geq0}$ is defined by
setting
\[
X_{k}= \zeta_{Y_{k}} \, .
\]
We assume that ${\mathbb{E}}[\tau]<\infty$ which ensures that $\{\phi
_{n}\}_{n\geq0}$ has a stationary distribution $\pi= (\pi_{i} , i \geq0)$
given by
\[
\pi_{j}=\pi_{0}\sum_{i=j+1}^{\infty}p_{i}\;,\;j=1,2\ldots
\]
where $\pi_{0}=1/{\mathbb{E}}(\tau)$. We also assume that $p_{j} >0$ for all
$j \geq0$. This last assumption implies the irreducibility of the Markov chain.

\medskip

The following notation will be  used: $ {\mathbb P}_{\phi_0=0}$ is the conditional probability given $\phi_0=0$ and $ {\mathbb E}_{\phi_0=0}$ is the expectation with respect to $ {\mathbb P}_{\phi_0=0}$.


In Corollary \ref{appliRWRTSstat} below, we shall make the following assumption
on the random time scenery:

\smallskip

\noindent\textbf{Condition ($A_{1}$)} \textit{$\{\zeta_{j}\}_{j\geq0}$ is a
strictly stationary sequence of centered random variables in ${\mathbb{L}}%
^{2}({\mathbb{P}})$, independent of $(\phi_{k})_{k\geq0}$ and such that
\begin{equation}
\sum_{k\geq1}\frac{\Vert{\mathbb{E}}(\zeta_{k}|\mathcal{G}_{0})\Vert_{2}%
}{\sqrt{k}}<\infty\,\text{ and }\lim_{n\rightarrow\infty}\sup_{j\geq i\geq
n}\Vert{\mathbb{E}}(\zeta_{i}\zeta_{j}|\mathcal{G}_{0})-{\mathbb{E}}(\zeta
_{i}\zeta_{j})\Vert_{1}=0\,,\label{condARWRTSpart1stat}%
\end{equation}
where }$\mathcal{G}$\textit{$_{i}=\sigma(\zeta_{k},k\leq i)$.}

\medskip

Applying Theorem \ref{CorollaryCLTNSforappli} and Proposition
\ref{toverifycondSqNS_g}, we can prove the following result concerning the
asymptotic behavior of $\{ n^{-1/2}S_{[nt]}, t \in[0,1] \} $ when the chain
starts from zero.

\begin{corollary}
\label{appliRWRTSstat} Assume that ${\mathbb{E}} (\tau^{2}) < \infty$ and that
$\{\zeta_{j}\}_{j \geq0}$ satisfies condition ($A_{1}$). Let $S_{0}=0$ and
$S_{k} = \sum_{i=1}^{k} X_{i}$ for any $k \geq1$. Then, under $ {\mathbb P}_{\phi_0=0}$, $\{ n^{-1/2}S_{[nt]}, t \in[0,1] \} $ converges in
distribution in $D[0,1]$ to a Brownian motion with parameter $c^{2}$ defined
by
\begin{equation}
\label{defc2stat}c^{2} = {\mathbb{E}} (\zeta_{0}^{2}) \Big ( 1 + 2 \sum_{i
\geq1} i \pi_{i} \Big ) + 2 \sum_{m \geq1} {\mathbb{E}} ( \zeta_{0} \zeta_{m}
) \sum_{j=1}^{m} (P^{j})_{0,m-j} \, .
\end{equation}

\end{corollary}

Note that ${\mathbb{E}} (\tau^{2}) < \infty$ is equivalent to $\sum_{i \geq1}
i \, {\mathbb{P}}(\tau> i) < \infty$ and therefore to the finitude of the
series $\sum_{i \geq1} i \pi_{i}$.


\section{Proofs} \label{sectionproofs}

In all the proofs, we use shall use the notation $a_{n} \ll b_{n}$ which means
that there exists a universal constant $C$ such that, for all $n \geq1$,
$a_{n} \leq C b_{n}$.

\subsection{Preparatory material}

The next result is a version of the functional central limit theorem for
triangular arrays of martingale differences essentially due to Aldous (1978)
and G\"{a}nssler-H\"{a}usler (1979) (see also Theorem 3.2 in Helland (1982)).

\begin{theorem}
[Aldous-G\"{a}nssler-H\"{a}usler]\label{martFCLTkn(t)} Let $v_{n}%
(\cdot):[0,1]\rightarrow\{0,1,2,\ldots,n\}$ be a sequence of integer valued,
non-decreasing and right-continuous functions. Assume $(d_{i,n})_{1\leq i\leq
n}$ is an array of martingale differences adapted to an array $(\mathcal{F}%
_{i,n})_{0\leq i\leq n}$ of nested sigma fields. Let $\sigma(\cdot)$ be a
non-negative function on $[0,1]$ such that $\sigma^{2}(\cdot)$ is Lebesgue
integrable. Suppose that the following conditions hold:
\begin{equation}
\max_{1\leq j\leq n}|d_{j,n}|\text{ is uniformly integrable,}%
\label{negl1fcltkn(t)}%
\end{equation}
and, for all $t\in\lbrack0,1]$,
\begin{equation}
\sum_{j=1}^{v_{n}(t)}d_{j,n}^{2}\rightarrow^{{\mathbb{P}}}\int_{0}^{t}%
\sigma^{2}(u)du\text{ }\ \text{as }n\rightarrow\infty
\,.\label{convcarrefcltkn(t)}%
\end{equation}
Then $\{\sum_{j=1}^{v_{n}(t)}d_{j,n},t\in\lbrack0,1]\}$ converges in
distribution in $D[0,1]$ to $\big \{\int_{0}^{t}\sigma(u)dW(u),t\in
\lbrack0,1]\big \}$ where $W$ is a standard Brownian motion.
\end{theorem}

\subsubsection{A maximal inequality in the non-stationary setting}

The following theorem is an extension of Proposition 2.3 in \cite{PU1} to the
non-stationary case. The proof follows the lines of the proof of Theorem 3 in
\cite{WuZh}, but in the non-stationary setting, and is then done by induction.
The proof is left to the reader but details can be found in the proof of
Theorem 3.2 in \cite{CDM}.

\begin{theorem}
\label{L2mainNS} Let $(X_{k})_{k\in{\mathbb{Z}}}$ be a sequence of real-valued
random variables in ${\mathbb{L}}^{2}$ and adapted to a filtration
$({\mathcal{F}}_{k})_{k\in{\mathbb{Z}}}$. Let $S_{n}=\sum_{k=1}^{n}X_{k}$,
$S_{0}=0$ and $S_{n}^{\ast}=\max_{1\leq k\leq n}|S_{k}|$. Then, for any
$n\geq1$,
\begin{equation}
\Vert S_{n}^{\ast}\Vert_{2}\leq3\Big (\sum_{j=1}^{n}\Vert X_{j}\Vert_{2}%
^{2}\Big )^{1/2}+3\sqrt{2}\Delta_{n}(X)\,, \label{inequality1maxSCpropCLTNS2}%
\end{equation}
where
\[
\Delta_{n}(X)=\sum_{j=0}^{r-1}\Big (\sum_{k=1}^{2^{r-j}}\Vert{\mathbb{E}%
}(S_{k2^{j}}-S_{(k-1)2^{j}}|{\mathcal{F}}_{(k-2)2^{j}+1})\Vert_{2}%
^{2}\Big )^{1/2}\, ,
\]
with $r$ the unique positive integer such that $2^{r-1}\leq n<2^{r}$.
\end{theorem}

\subsection{Proof of Theorem \ref{fclt_MPU}}

Recall that $X:=\{X_{k,n}:k=1,\ldots,n\}=(X_{k,n})_{k=1}^{n}$ is a triangular
array of real-valued random variables in ${\mathbb{L}}^{2}$ adapted to a
filtration $(\mathcal{F}_{k,n})_{0\leq k\leq n}$. Without loss of generality,
we assume that $X_{k,n}=0$ for $k>n$ and $\mathcal{F}_{k,n}=\mathcal{F}_{n,n}$
for $k>n$. Moreover, by abuse of notation, we will often avoid the index $n$.
In particular, we shall write $X_{k}=X_{k,n}$ and $\mathcal{F}_{k}%
=\mathcal{F}_{k,n}$, and we will use the notations
\[
{\mathbb{E}}_{j}(X)={\mathbb{E}}(X|\mathcal{F}_{j})\,,\,\mathbf{P}%
_{j}(x)={\mathbb{E}}_{j}(X)-{\mathbb{E}}_{j-1}(X)\,.
\]
Fix positive integer $n$ and define the unique positive integer $r$ such that
$2^{r-1}\leq n<2^{r}$.

For each $n$, we define the following random variables:
\[
S_{n}=S_{n}(X)=\sum_{k=1}^{n}X_{k}\,,\,S_{0}=0\,,\,S_{n}^{\ast}=S_{n}^{\ast
}(X)=\max_{1\leq k\leq n}|S_{k}|\,.\
\]

Theorem \ref{fclt_MPU} will follow from a martingale approximation and an
application of Theorem \ref{martFCLTkn(t)}, for the approximating martingale.

\subsubsection{Step 1: A general Lemma.}

Let us first introduce some notations. Let $m$ be a fixed positive integer
such that $m<n$. Let us then define
\begin{equation}
\theta_{{\ell}}^{m}=\frac{1}{m}\sum_{i=1}^{m-1}{\mathbb{E}}_{{\ell}}%
(X_{{\ell+1}}+\ldots+X_{{\ell}+i})\,,\,D_{{\ell}}^{m}=\frac{1}{m}\sum_{i={0}%
}^{m-1}{\mathbf{P}}_{\ell}(S_{{\ell}+i})=\frac{1}{m}\sum_{i={0}}%
^{m-1}{\mathbf{P}}_{\ell}(S_{{\ell}+i}-S_{\ell-1}) \,, \label{martdefNS_1}%
\end{equation}
and
\begin{equation}
Y_{\ell}^{m}=\frac{1}{m}{\mathbb{E}}_{{\ell}}(S_{{\ell}+m}-S_{{\ell}%
})\,,\,R_{k}^{m}=\sum_{{\ell}=0}^{k-1}Y_{\ell}^{m}\,. \label{defYNS_1}%
\end{equation}
Then, $D^{m}=(D_{k}^{m})_{k=1}^{n}$ is a (triangular) array of martingale
differences adapted to the filtration $(\mathcal{F}_{{k}})_{0\leq k\leq n}$
and the following decomposition is valid:
\begin{equation}
X_{{\ell}}=D_{{\ell}}^{m}+\theta_{{\ell-1}}^{m}-\theta_{{\ell}}^{m}+Y_{\ell
-1}^{m}\,. \label{martrepresXkNS_1}%
\end{equation}
Also, for any positive integer $m$ and $k$, we have
\begin{equation}
S_{k}=M_{k}^{m}+\theta_{0}^{m}-\theta_{k}^{m}+R_{k}^{m}\,. \label{martdecNS_1}%
\end{equation}
As an intermediate step in proving Theorem \ref{fclt_MPU} we shall prove a
lemma under a set of assumptions which will be verified later. The next
assumption \textbf{$(H)$} aims to guarantee that, in a certain sense, $S_{k}$
can be approximated by $M_{k}^{m^{\prime}}$ (for $m^{\prime}$ a subsequence of
$m$) and it is then used to verify the conditions of Theorem
\ref{martFCLTkn(t)}.

There exists an increasing subsequence of integers $(m_{j})_{j\geq1}$ with
$m_{j}\rightarrow\infty$ as $j\rightarrow\infty$ such that
\[
(H):=%
\begin{cases}
\lim_{j\rightarrow\infty}\sup_{n\geq1}\sum_{\ell=0}^{n-1}\Vert Y_{\ell}%
^{m_{j}}\Vert_{2}^{2}=0\,,\\
\lim_{j\rightarrow\infty}\sup_{n\geq1}\Delta_{n}(Y^{m_{j}})=0\,,\\
\lim_{j\rightarrow\infty}\sup_{n\geq1}\sum_{k=0}^{n-1}\Vert\theta_{k}^{m_{j}%
}\Vert_{2}\Vert Y_{k}^{m_{j}}\Vert_{2}=0\,,
\end{cases}
\]
where
\[
\Delta_{n}(Y^{m}):=\sum_{r=0}^{d}\Big (\sum_{k=1}^{2^{d-r}}\Vert{\mathbb{E}%
}(R_{k2^{r}}^{m}-R_{(k-1)2^{r}}^{m}|{\mathcal{F}}_{(k-2)2^{r}})\Vert_{2}%
^{2}\Big )^{1/2}\,.
\]
We are now in the position to state our general lemma.

\begin{lma}
\label{Gen_fclt} Assume that the Lindeberg-type condition \eqref{condLNS_g}
holds and that condition $(H)$ is satisfied. Assume in addition that there
exist a sequence of non-decreasing and right-continuous functions $v_{n}%
(\cdot):[0,1]\rightarrow\{1,2,\ldots,n\}$ and a non-negative real $c^{2}$ such
that condition \eqref{AssGr3_1_g} holds. Then $\big \{\sum_{k=1}^{v_{n}%
(t)}X_{k},t\in\lbrack0,1]\big \}$ converges in distribution in $D([0,1])$ to
$cW$ where $W$ is a standard Brownian motion.
\end{lma}

\noindent\textbf{Proof.} To soothe the notations, we will often write $m$
instead of $m_{j}$. To prove the lemma, let us first analyze the negligibility
in some sense of the variables $\theta_{k}^{m}$ and $R_{k}^{m}$. Notice that
from the definition \eqref{martdefNS_1}
\begin{align*}
\max_{0\leq k\leq n}|\theta_{k}^{m}|^{2} &  =\max_{0\leq k\leq n}\Big(\frac
{1}{m}\sum_{i=1}^{m-1}{\mathbb{E}}_{k}(X_{k+1}+\ldots+X_{k+i})\Big)^{2}\\
&  \leq m^{2}\max_{0\leq j\leq n}\mathbb{E}_{j}\big(\max_{1\leq k\leq n}%
|X_{k}|^{2}\big).
\end{align*}
By applying the Doob's maximal inequality and next truncation, we derive
\begin{align*}
\mathbb{E}\big[\max_{0\leq j\leq n}\mathbb{E}_{j}\big(\max_{1\leq k\leq
n}|X_{k}|^{2}\big)\big] &  \leq4\mathbb{E}\Big (\mathbb{E}_{n}\big(\max_{1\leq
k\leq n}|X_{k}|^{2}\big)\Big )=4\mathbb{E}\big(\max_{1\leq k\leq n}|X_{k}%
|^{2}\big)\\
&  \leq4\varepsilon+4\sum_{k=1}^{n}{\mathbb{E}}\{X_{k}^{2}I(|X_{k}%
|>\varepsilon)\}.
\end{align*}
Combining it with the previous estimate, taking into account the
Lindeberg-like condition \eqref{condLNS_g} and letting $n$ tend to infinity
and then $\varepsilon\rightarrow0$ we obtain for each $m$, that
\begin{equation}
\mathbb{E}\big(\max_{0\leq k\leq n}|\theta_{k}^{m}|^{2}\big)\rightarrow
0\mbox{ as }n\rightarrow\infty\,.\label{theta_1}%
\end{equation}
Note that, proceeding similarly, we also have that, for each $m$,
\begin{equation}
\mathbb{E}\big(\max_{1\leq k\leq n}|D_{k}^{m}|^{2}\big)\rightarrow
0\mbox{ as }n\rightarrow\infty\,.\label{Diff_1}%
\end{equation}

Now, by applying Theorem \ref{L2mainNS} to the array $(Y_{k}^{m}%
)_{k\in{\mathbb{Z}}}$, we have
\[
\big \Vert\max_{1\leq k\leq n}|R_{k}^{m}|\big \Vert_{2}=\big \Vert\max_{1\leq
k\leq n}\big |\sum_{\ell=0}^{k-1}Y_{\ell}^{m}\big |\big \Vert_{2}%
\leq3\Big (\sum_{k=0}^{n-1}\Vert Y_{k}^{m}\Vert_{2}^{2}\Big )^{1/2}+3\sqrt
{2}\Delta_{n}(Y^{m})\,.
\]
Taking $m=m_{j}$, by assumption $(H)$ the terms in the r.h.s tend to $0,$
uniformly in $n$ by letting $j\rightarrow\infty$. Hence, we derive the bound
\begin{equation}
\sup_{n\geq1}\big \Vert\max_{1\leq k\leq n}|R_{k}^{m_{j}}|\big \Vert_{2}%
\rightarrow0\mbox{ as }j\rightarrow\infty\,.\label{max_R_m_1}%
\end{equation}
By the relations (\ref{theta_1}) and (\ref{max_R_m_1}) we have the following
martingale approximation%
\[
\limsup_{n}\Big \Vert \max_{1\leq k\leq n} \Big |\sum_{i=1}^{k}X_{i,n}-\sum_{\ell=1}%
^{k}D_{{\ell}}^{m_{j}} \Big |\Big \Vert_{2}\rightarrow0\text{ as }j\rightarrow
\infty\,.
\]
This limit clearly implies%
\begin{equation}
\limsup_{n}\Big \Vert \sup_{t\in\lbrack0,1]} \Big |\sum_{i=1}^{v_{n}(t)}X_{i,n}%
-\sum_{\ell=1}^{v_{n}(t)}D_{{\ell}}^{m_{j}} \Big |\Big \Vert_{2}\rightarrow0\text{ as
}j\rightarrow\infty\,,\label{mart approx}%
\end{equation}
and also for $j_{0}$ fixed, 
\[
\sup_{n\geq1}{\mathbb{E}}(S_{n}^{2})\leq\sup_{n\geq1} \sum_{\ell=1}%
^{n}\big \Vert D_{{\ell}}^{m_{j_{0}}}\big \Vert_{2} + \varepsilon_{j_0}\,,
\]
where $\varepsilon_{j_0}$ is a finite positive constant. 
Now, by definition \eqref{martdefNS_1},
\[
\Vert D_{k}^{m_{j_{0}}}\Vert_{2}\leq\frac{1}{m_{j_{0}}}\sum_{i={0}}^{m_{j_{0}%
}-1}\Vert{\mathbf{P}}_{k}(S_{k+i}-S_{k})\Vert_{2}\leq\frac{1}{m_{j_{0}}}%
\sum_{i={0}}^{m_{j_{0}}-1}\Vert{\mathbb{E}}_{k}(S_{k+i}-S_{k})\Vert_{2}\,.
\]
Hence, since $X_{k}=X_{k,n}=0$, $k>n$,
\begin{equation}
\sum_{k=1}^{n}\Vert D_{k}^{m_{j_{0}}}\Vert_{2}^{2}\leq m_{j_{0}}^{2}\sum
_{k=1}^{n}\Vert X_{k}\Vert_{2}^{2}\,.\label{DLNS_g}%
\end{equation}
Therefore, by the first part of \eqref{condLNS_g}, 
\begin{equation}
\sup_{n\geq1}{\mathbb{E}}(S_{n}^{2})\leq C_{j_{0}}<\infty
\,.\label{FCLTvariance}%
\end{equation}

From (\ref{theta_1}), \eqref{Diff_1}, (\ref{max_R_m_1}), (\ref{mart approx})
and (\ref{AssGr3_1_g}), via Lemma \ref{lmaconvergence} in the Appendix, we
deduce that we can find a sequence of positive integers $\ell(n)$ such that
$\ell(n)\rightarrow\infty$ and setting $m_{n}^{\prime}=m_{\ell(n)}$,
\begin{equation}
\lim_{n\rightarrow\infty}\big \Vert\max_{0\leq k\leq n}|\theta_{k}%
^{m_{n}^{\prime}}|\big \Vert_{2}=0\,,\,\lim_{n\rightarrow\infty}%
\big \Vert\max_{1\leq k\leq n}|D_{k}^{m_{n}^{\prime}}|\big \Vert_{2}%
=0\,,\label{theta_1(2)}%
\end{equation}

\begin{equation}
\lim_{n\rightarrow\infty}\big \Vert\max_{1\leq k\leq n}|R_{k}^{m_{n}^{\prime}%
}|\big \Vert_{2}=0 \, ,\label{max_R_m_1(2)}%
\end{equation}

\begin{equation}
\lim_{n\rightarrow\infty} \Big \Vert \sup_{t\in\lbrack0,1]} \Big |\sum_{i=1}^{v_{n}(t)} 
X_{i,n}-\sum_{\ell=1}^{v_{n}(t)}D_{{\ell}}^{m_{n}^{\prime}} \Big |
\Big \Vert_{2}=0\, , \label{mart approx(2)}%
\end{equation}
and, for any $t \in [0,1]$, 
\begin{equation}
\ \lim_{n\rightarrow\infty}{\mathbb{P}}\Big (\Big |\sum_{k=1}^{v_{n}%
(t)}\big (X_{k}^{2}+2X_{k}\theta_{k}^{m_{n}^{\prime}}\big )-tc^{2}%
\Big |>\varepsilon\Big )=0\,. \label{AssGr3_1_g(2)}%
\end{equation}

In addition, by condition $(H)$, on the same subsequence $(m_{n}^{\prime})$ we
also have%

\[
(H^{\prime}):=%
\begin{cases}
\lim_{n\rightarrow\infty}\sum_{\ell=0}^{n-1}\Vert Y_{\ell}^{m_{n}^{\prime}%
}\Vert_{2}^{2}=0\,,\\
\lim_{n\rightarrow\infty}\Delta_{n}(Y^{m_{n}^{\prime}})=0\,,\\
\lim_{n\rightarrow\infty}\sum_{k=0}^{n-1}\Vert\theta_{k}^{m_{n}^{\prime}}%
\Vert_{2}\Vert Y_{k}^{m_{n}^{\prime}}\Vert_{2}=0\,.
\end{cases}
\]

By \eqref{mart approx(2)}, it suffices to show that  $\big \{ \sum_{\ell=1}^{v_{n}(t)}D_{{\ell}}^{m_{n}^{\prime}} ,t\in\lbrack0,1]\big \}$ converges in distribution in $D([0,1])$ to
$cW$. We shall verify now that the triangular array of martingale differences $(D_{{\ell}}^{m_{n}^{\prime}})_{1 \leq \ell\leq n }$ satisfies the conditions of Theorem \ref{martFCLTkn(t)}. The
condition (\ref{negl1fcltkn(t)}) follows from the second part of
(\ref{theta_1(2)}). In order to verify condition (\ref{convcarrefcltkn(t)}%
)\ we proceed in the following way. We start from the identity
(\ref{martrepresXkNS_1}) written as ($m=m_{n}^{\prime}$)%
\[
X_{{\ell}}+\theta_{{\ell}}^{m}=D_{{\ell}}^{m}+\theta_{{\ell-1}}^{m}+Y_{\ell
-1}^{m}\,.
\]
Therefore%
\[
X_{{\ell}}^{2}+2X_{{\ell}}\theta_{{\ell}}^{m}+(\theta_{{\ell}}^{m}%
)^{2}=(D_{{\ell}}^{m})^{2}+(\theta_{{\ell-1}}^{m})^{2}+(Y_{\ell-1}^{m}%
\,)^{2}+2\theta_{{\ell-1}}^{m}Y_{\ell-1}^{m}+2D_{{\ell}}^{m}(\theta_{{\ell-1}%
}^{m}+Y_{\ell-1}^{m})\,.
\]
We sum over $\ell$ and get
\[
\sum\nolimits_{\ell=1}^{v_{n}(t)}(X_{{\ell}}^{2}+2X_{{\ell}}\theta_{{\ell}%
}^{m})+(\theta_{v_{n}(t)}^{m})^{2}=\sum\nolimits_{\ell=1}^{v_{n}(t)}(D_{{\ell
}}^{m})^{2}+(\theta_{{0}}^{m})^{2}+\sum\nolimits_{\ell=1}^{v_{n}(t)}2D_{{\ell
}}^{m}(\theta_{{\ell-1}}^{m}+Y_{\ell-1}^{m})+R^{\prime}(v_{n}(t))\,,
\]
where%
\[
R^{\prime}(v_{n}(t))=\sum_{\ell=0}^{v_{n}(t)-1}(Y_{\ell}^{m_{n}^{\prime}}%
)^{2}+2\sum_{k=0}^{v_{n}(t)-1}\theta_{k}^{m_{n}^{\prime}}Y_{k}^{m_{n}^{\prime
}}\,.
\]
By the Cauchy-Schwarz inequality and condition $(H^{\prime})$ we have that%
\[
\mathbb{E}(\sup_{0\leq t\leq1}|R^{\prime}(v_{n}(t))|)\leq\sum_{\ell=0}%
^{n-1}\Vert Y_{\ell}^{m_{n}^{\prime}}\Vert_{2}^{2}+2\sum_{k=0}^{n-1}%
\Vert\theta_{k}^{m_{n}^{\prime}}\Vert_{2}\Vert Y_{k}^{m_{n}^{\prime}}\Vert
_{2}\rightarrow0\text{ as }n\rightarrow\infty\,.
\]
Furthermore, by using the first part of (\ref{theta_1(2)}) we also have
\[
\mathbb{E}\sup_{0\leq t\leq1}|(\theta_{v_{n}(t)}^{m})^{2}-(\theta_{{0}}%
^{m})^{2}|\rightarrow0\text{ as }n\rightarrow\infty\,.
\]
Now, by gathering the above considerations and by also using
(\ref{AssGr3_1_g(2)}), we shall have
\[
\sum\nolimits_{\ell=1}^{v_{n}(t)}(D_{{\ell}}^{m_{n}^{\prime}})^{2}\rightarrow
c^{2}t\text{ in probability as }n\rightarrow\infty\,,
\]
if we can prove that $\sum\nolimits_{\ell=1}^{v_{n}(t)}2D_{{\ell}}%
^{m_{n}^{\prime}}(\theta_{{\ell-1}}^{m_{n}^{\prime}}+Y_{\ell-1}^{m_{n}%
^{\prime}})\rightarrow0$ in probability. Because $(\theta_{{\ell-1}}%
^{m_{n}^{\prime}}+Y_{\ell-1}^{m_{n}^{\prime}})$ is a previsible (i.e.
$\mathcal{F}_{\ell-1,n}$-measurable) random variable, the result follows again
from (\ref{theta_1(2)}) and $(H^{\prime})$, by using the following fact, which
is Theorem 2.11 in Hall-Heyde (1980):

\begin{fact}
\label{LLN_trian} Let $(Z_{i})_{i=1}^{n}$ be real-valued martingale
differences adapted to a non-increasing filtration $(\mathcal{F}_{i})_{0\leq
i\leq n}$ and let $(A_{k})_{k=1}^{n}$ be real-valued random variables such
that $A_{k}$ is $\mathcal{F}_{k-1}$-measurable. Then, there exists a positive
constant $c$ such that
\[
{\mathbb{E}}\max_{1\leq k\leq n}\Big|\sum_{i=1}^{k}A_{i}Z_{i}\Big|\leq
c\big \{{\mathbb{E}}\max_{1\leq k\leq n}|A_{k}|^{2}\big \}^{1/2}%
\Big \{\sum_{i=1}^{n}{\mathbb{E}}(Z_{i}^{2})\Big \}^{1/2}\,.
\]

\end{fact}

together with the following remark: by \eqref{mart approx(2)} and
\eqref{FCLTvariance},
\begin{multline*}
\limsup_{n}\sum_{\ell=1}^{v_{n}(t)}\Vert D_{{\ell}}^{m_{n}^{\prime}}\Vert
_{2}^{2}=\limsup_{n}\Big \Vert\sum_{\ell=1}^{v_{n}(t)}D_{{\ell}}%
^{m_{n}^{\prime}}\Big \Vert_{2}^{2}\\
\leq\limsup_{n\rightarrow\infty}\sup_{t\in\lbrack0,1]}\big \Vert\sum
_{i=1}^{v_{n}(t)}X_{i,n}-\sum_{i=1}^{v_{n}(t)}D_{{\ell}}^{m_{n}^{\prime}%
}\big \Vert_{2}^{2}+\limsup_{n\rightarrow\infty}\big \Vert\sum_{i=1}%
^{v_{n}(t)}X_{i,n}\big \Vert_{2}^{2}\leq C_{j_{0}}\,.
\end{multline*}
This ends the proof of the lemma. $\diamond$

\subsubsection{Step 2: end of the proof of Theorem \ref{fclt_MPU}.}

We are going to prove that Theorem \ref{fclt_MPU} follows from an application
of Lemma \ref{Gen_fclt}. With this aim we start by noticing the following
fact: if the second part of \eqref{MWcondNS2_g} holds then there exists an
increasing subsequence of integers $(m(j))_{j \geq1}$ with $m(j)
\rightarrow\infty$ as $j \rightarrow\infty$ and such that
\begin{equation}
\label{consequenceMWcondNS2}\lim_{j \rightarrow\infty} \sum_{\ \ell\geq m(j)}
2^{- \ell/2} B ( 2^{\ell}, 2^{m(j)} ) =0 \, .
\end{equation}

Hence, to show that condition $(H)$ of Lemma \ref{Gen_fclt} holds, we shall
prove that its three assumptions are satisfied with $m_{j}=2^{m(j)}$. So, in
what follows $m_{j}=2^{m(j)}$ where $m(j)$ is an increasing subsequence of
integers tending to infinity and such that \eqref{consequenceMWcondNS2} holds.
As before, we will sometimes write $m$ instead of $m_{j}$.

\bigskip

\textbf{Veryfying first condition in $(H)$}. We first notice that, by the
definition \eqref{defalu_g} and first part of condition \eqref{MWcondNS2_g}
\begin{equation}
\sup_{n\geq1}\sum_{k=0}^{n-1}\Vert Y_{k}^{m_{j}}\Vert_{2}^{2}=m_{j}^{-2}%
\sup_{n\geq1}\sum_{k=0}^{n-1}\Vert{\mathbb{E}}_{k}(S_{k+m_{j}%
}-S_{k})\Vert_{2}^{2} = m_{j}^{-2}A^{2}(m_{j})\,\rightarrow
0\mbox{ as }j\rightarrow\infty\,, \label{p1plusNS_g}%
\end{equation}
which proves the first condition in $(H)$.

\medskip\textbf{Veryfying second condition in $(H)$}. This needs more
considerations. It is convenient to use the decomposition
\begin{multline*}
\sum_{r=0}^{d}\Big (\sum_{k=1}^{2^{d-r}}\Vert{\mathbb{E}}(R_{k2^{r}}%
^{m}-R_{(k-1)2^{r}}^{m}|{\mathcal{F}}_{(k-2)2^{r}})\Vert_{2}^{2}%
\Big )^{1/2}=\sum_{r=0}^{b}\Big (\sum_{k=1}^{2^{d-r}}\Vert{\mathbb{E}%
}(R_{k2^{r}}^{m}-R_{(k-1)2^{r}}^{m}|{\mathcal{F}}_{(k-2)2^{r}})\Vert_{2}%
^{2}\Big )^{1/2}\\
+\sum_{r=b+1}^{d}\Big (\sum_{k=1}^{2^{d-r}}\Vert{\mathbb{E}}(R_{k2^{r}}%
^{m}-R_{(k-1)2^{r}}^{m}|{\mathcal{F}}_{(k-2)2^{r}})\Vert_{2}^{2}%
\Big )^{1/2}\,,
\end{multline*}
where $b$ is the unique positive integer such that $2^{b}\leq m<2^{b+1}$. To
estimate the first sum in the right-hand side, notice that, by the properties
of the conditional expectation, we have
\begin{align}
\Vert{\mathbb{E}}(R_{k2^{r}}^{m}-R_{(k-1)2^{r}}^{m}  &  |{\mathcal{F}%
}_{(k-2)2^{r}})\Vert_{2}\leq\sum_{{\ell}=0}^{2^{r}-1}\Vert{\mathbb{E}}%
(Y_{\ell+(k-1)2^{r}}^{m}|{\mathcal{F}}_{(k-2)2^{r}})\Vert_{2}%
\nonumber\label{dec1R30-10_g}\\
&  \leq\frac{1}{m}\sum_{{\ell}=0}^{2^{r}-1}\Vert{\mathbb{E}}(S_{\ell
+(k-1)2^{r}+m}-S_{\ell+(k-1)2^{r}}|{\mathcal{F}}_{(k-2)2^{r}})\Vert
_{2}\nonumber\\
&  \leq\frac{1}{m}\sum_{{\ell}=0}^{2^{r}-1}\Vert{\mathbb{E}}(S_{\ell
+(k-1)2^{r}+m}-S_{\ell+(k-1)2^{r}}|{\mathcal{F}}_{(k-2)2^{r}+\ell})\Vert
_{2}\nonumber\\
&  \leq\frac{1}{m}\sum_{{\ell}=(k-1)2^{r}}^{k2^{r}-1}\Vert{\mathbb{E}}%
(S_{\ell+m}-S_{\ell}|{\mathcal{F}}_{\ell-2^{r}})\Vert_{2}\,.
\end{align}
Therefore, by definition \eqref{defalu_g},
\begin{multline*}
\Big (\sum_{k=1}^{2^{d-r}}\Vert{\mathbb{E}}(R_{k2^{r}}^{m}-R_{(k-1)2^{r}}%
^{m}|{\mathcal{F}}_{(k-2)2^{r}})\Vert_{2}^{2}\Big )^{1/2}\leq\frac{2^{r/2}}%
{m}\Big (\sum_{k=1}^{2^{d-r}}\sum_{{\ell}=(k-1)2^{r}}^{k2^{r}-1}%
\Vert{\mathbb{E}}(S_{\ell+m}-S_{\ell}|{\mathcal{F}}_{\ell-2^{r}})\Vert_{2}%
^{2}\Big )^{1/2}\\
\leq\frac{2^{r/2}}{m}\Big (\sum_{{\ell}=0}^{2^{d}-1}\Vert{\mathbb{E}}%
(S_{\ell+m}-S_{\ell}|{\mathcal{F}}_{\ell-2^{r}})\Vert_{2}^{2}\Big )^{1/2}%
\leq\frac{2^{r/2}}{m}A(m)\,,
\end{multline*}
giving
\[
\sum_{r=0}^{b}\Big (\sum_{k=1}^{2^{d-r}}\Vert{\mathbb{E}}(R_{k2^{r}}%
^{m}-R_{(k-1)2^{r}}^{m}|{\mathcal{F}}_{(k-2)2^{r}})\Vert_{2}^{2}%
\Big )^{1/2}\leq\frac{2}{\sqrt{2}-1}\frac{A(m)}{\sqrt{m}}\,.
\]
To estimate the second sum we also apply the properties of the conditional
expectation and write this time
\begin{align}
\Vert{\mathbb{E}}(R_{k2^{r}}^{m}  &  -R_{(k-1)2^{r}}^{m}|{\mathcal{F}%
}_{(k-2)2^{r}})\Vert_{2}\nonumber\label{dec2R30-10_g}\\
&  \leq\Big \Vert\frac{1}{m}\sum_{u=0}^{m-1}{\mathbb{E}}(S_{k2^{r}%
+u}-S_{(k-1)2^{r}+u}|{\mathcal{F}}_{(k-2)2^{r}})\Big \Vert_{2}\\
&  \leq\Big \Vert\frac{1}{m}\sum_{u=0}^{m-1}{\mathbb{E}}(S_{k2^{r}%
+u}-S_{(k-1)2^{r}+u}|{\mathcal{F}}_{(k-2)2^{r}+1})\Big \Vert_{2}:=\Vert
{\bar{S}}_{k-1}(2^{r},m)\Vert_{2}\,.\nonumber
\end{align}
Hence, by definition \eqref{defblu_g}
\[
\sum_{r=b+1}^{d}\Big (\sum_{k=1}^{2^{d-r}}\Vert{\mathbb{E}}(R_{k2^{r}}%
^{m}-R_{(k-1)2^{r}}^{m}|{\mathcal{F}}_{(k-2)2^{r}})\Vert_{2}^{2}%
\Big )^{1/2}\leq\sum_{r=b+1}^{d}B(2^{r},m)\,.
\]
So, overall,
\begin{equation}
\sum_{r=0}^{d}\Big (\sum_{k=1}^{2^{d-r}}\Vert{\mathbb{E}}(R_{k2^{r}}%
^{m}-R_{(k-1)2^{r}}^{m}|{\mathcal{F}}_{(k-2)2^{r}})\Vert_{2}^{2}%
\Big )^{1/2}\leq\frac{2}{\sqrt{2}-1}\frac{A(m)}{\sqrt{m}}+\sqrt{2}\sum
_{r=b+1}^{d}B(2^{r},m)\,. \label{p2plusNS_g}%
\end{equation}
This gives
\[
\sup_{n\geq1}\Delta_{n}(Y^{m_{j}})\leq\frac{2}{\sqrt{2}-1}2^{-m(j)/2}%
A(2^{m(j)})+\sqrt{2}\sum_{r\geq m(j)}B(2^{r},2^{m(j)})
\]
which, together with condition \eqref{MWcondNS2_g}, prove the second condition
in $(H)$. \medskip

It is worth to notice that we have proved the following maximal inequality for
the array of $Y$'s.

\begin{lemma}
\label{plusNS_g}There exists a positive constant $C$ such that, for every
positive integers $n$ and $m$ such that $m \leq n$,
\[
\big \Vert\max_{1\leq j\leq n} \big |\sum_{k=0}^{j-1}Y_{k}^{m}
\big | \big \Vert_{2} \leq3 \Big ( 1 + \frac{ 2 \sqrt{2} }{\sqrt2 -1}
\Big ) \frac{A(m)}{ \sqrt m }+ 6 \sum_{r=[\log_{2} (m) ]}^{d} 2^{-r/2}
b(2^{r},m) \, ,\label{normNS}%
\]
where $d$ be the unique positive integer such that $2^{d-1} \leq n < 2^{d}$.
\end{lemma}

\noindent\textit{Proof of Lemma \ref{plusNS_g}.} We start with the general
maximal inequality \eqref{inequality1maxSCpropCLTNS2} applied to the array
$(Y^{m}_{k})_{k \in{\mathbb{Z}}}$. Let $d$ be the unique positive integer such
that $2^{d-1} \leq n < 2^{d}$ (so $d \leq[ \log_{2} (2n) ]$). Recalling the
notation $R_{k}^{m}=\sum_{{\ell}=0}^{k-1}Y_{\ell}^{m} $, we get
\begin{equation}
\label{inequality1maxSCpropCLTNS2forY}\big \Vert\max_{1\leq j\leq n}
\big |\sum_{k=0}^{j-1}Y_{k}^{m} \big | \big \Vert_{2} \leq3 \Big ( \sum
_{k=0}^{n-1} \Vert Y_{k}^{m} \Vert_{2}^{2} \Big )^{1/2} + 3 \sqrt2 \sum
_{r=0}^{d} \Big ( \sum_{k=1}^{2^{d-r}} \Vert{\mathbb{E}} ( R^{m}_{k 2^{r} }
-R^{m}_{(k-1) 2^{r} } | {\mathcal{F}}_{(k-2) 2^{r} } ) \Vert_{2}^{2}
\Big )^{1/2} \, .
\end{equation}
It remains to apply bounds \eqref{p1plusNS_g} and \eqref{p2plusNS_g}.
$\diamond$

\medskip\textbf{Veryfying third condition in $(H)$}. For any positive integer
$i$ such that $i<m_{j}$, we write its decomposition in basis $2$,
\[
i=\sum_{k=0}^{[\log_{2}(i)]+1}c_{k}(i)2^{k}\ \text{ where $c_{k}(i)\in
\{0,1\}$}\,.
\]
Denote by $i_{u}=\sum_{k=0}^{u}c_{k}(i)2^{k}$ (hence $i_{[\log_{2}(i)]+1}=i$),
for $u\geq0$ and set $i_{-1}=0$. We have
\begin{multline*}
\Vert\theta_{{\ell}}^{m}\Vert_{2}\leq\frac{1}{m}\sum_{i=1}^{m-1}%
\Vert{\mathbb{E}}_{\ell}(S_{\ell+i}-S_{\ell})\Vert_{2}=\frac{1}{m}\sum
_{i=1}^{m-1}\Big \Vert\sum_{u=0}^{[\log_{2}(i)]+1}{\mathbb{E}}_{\ell}%
(S_{\ell+i_{u}}-S_{\ell+i_{u-1}})\Big \Vert_{2}\\
\leq\frac{1}{m}\sum_{i=1}^{m-1}\sum_{u=0}^{[\log_{2}(i)]+1}\Vert{\mathbb{E}%
}_{\ell+i_{u-1}}(S_{\ell+i_{u}}-S_{\ell+i_{u-1}})\Vert_{2}\\
=\frac{1}{m}\sum_{i=1}^{m-1}\sum_{u=0}^{[\log_{2}(i)]+1}c_{u}(i)\Vert
{\mathbb{E}}_{\ell+i_{u-1}}(S_{\ell+i_{u-1}+2^{u}}-S_{\ell+i_{u-1}})\Vert
_{2}\,.
\end{multline*}
Hence, by taking into account definition \eqref{defalu_g},
\begin{multline}
\sum_{\ell=0}^{n-1}\Vert\theta_{{\ell}}^{m}\Vert_{2}\Vert Y_{\ell}^{m}%
\Vert_{2}\label{boundviabasis2_g}\\
\leq\frac{1}{m}\sum_{i=1}^{m-1}\sum_{u=0}^{[\log_{2}(i)]+1}\Big (\sum_{\ell
=0}^{n-1}c_{u}(i)\Vert{\mathbb{E}}_{\ell+i_{u-1}}(S_{\ell+i_{u-1}+2^{u}%
}-S_{\ell+i_{u-1}})\Vert_{2}^{2}\Big )^{1/2}\Big (\sum_{\ell=0}^{n-1}\Vert
Y_{\ell}^{m}\Vert_{2}^{2}\Big )^{1/2}\\
\leq\frac{1}{m}\sum_{i=1}^{m-1}\sum_{u=0}^{[\log_{2}(i)]+1}A(2^{u}%
)\Big (\sum_{\ell=0}^{n-1}\Vert Y_{\ell}^{m}\Vert_{2}^{2}\Big )^{1/2}\,.
\end{multline}
So, by the first part of condition \eqref{MWcondNS2_g}, there exists a
constant $C$ such that
\begin{multline*}
\sum_{\ell=0}^{n-1}\Vert\theta_{{\ell}}^{m}\Vert_{2}\Vert Y_{\ell}^{m}%
\Vert_{2}\leq\frac{C}{m}\sum_{i=1}^{m-1}\sum_{u=0}^{[\log_{2}(i)]+1}%
2^{u/2}\Big (\sum_{\ell=0}^{n-1}\Vert Y_{\ell}^{m}\Vert_{2}^{2}\Big )^{1/2}\\
\leq\frac{2C}{\sqrt{2}-1}\frac{1}{m}\sum_{i=1}^{m-1}\sqrt{i}\Big (\sum
_{\ell=0}^{n-1}\Vert Y_{\ell}^{m}\Vert_{2}^{2}\Big )^{1/2}\leq\frac{2C}%
{\sqrt{2}-1}\sqrt{m}\Big(\sum_{\ell=0}^{n-1}\Vert Y_{\ell}^{m}\Vert_{2}%
^{2}\Big )^{1/2}\,.
\end{multline*}
With $m=m_{j}=2^{m(j)}$ and taking now into account \eqref{p1plusNS_g}, it
follows that
\[
\sum_{\ell=0}^{n-1}\Vert\theta_{{\ell}}^{m_{j}}\Vert_{2}\Vert Y_{\ell}^{m_{j}%
}\Vert_{2}\leq\frac{2C}{\sqrt{2}-1}\frac{1}{\sqrt{m_{j}}}A(m_{j})=\frac
{2C}{\sqrt{2}-1}2^{-m(j)/2}A(2^{m(j)})\,,
\]
which converges to zero as $j\rightarrow\infty$ by the first part of
\eqref{MWcondNS2_g}. This shows that the third condition in $(H)$ is satisfied
and ends the proof of the theorem. $\diamond$

\subsection{Proof of Proposition \ref{toverifycondSqNS_g}}

Once again, to soothe the notation, we will avoid the index $n$ involved in
the variables and in the $\sigma$-algebras. In particular, we shall write
$X_{k}=X_{k,n}$ and $\mathcal{F}_{k}=\mathcal{F}_{k,n}$, and we will use the
notations ${\mathbb{E}}_{j} (X) = {\mathbb{E}}(X|\mathcal{F}_{j}) $ and
$\mathbf{P}_{j} (x) = {\mathbb{E}}_{j} (X) - {\mathbb{E}}_{j-1} (X)$. Moreover, without loss of generality,
we assume that $X_{k,n}=0$ for $k>n$

\smallskip

Clearly it is enough to show that, for any $t \in[0,1]$ and any fixed
non-negative integer $\ell$,
\begin{equation}
\label{cequonveutpourlesP0butFlo}\lim_{n \rightarrow\infty} \Big \Vert \sum
_{k=1}^{v_{n}(t)} \big ( {\mathbb{E}}_{k} ( X_{k}X_{k+\ell} ) - {\mathbb{E}%
}_{0} ( X_{k}X_{k+\ell} ) \big ) \Big \Vert_{1} = 0 \, .
\end{equation}
With this aim, we first write the following decomposition: for any positive
fixed integer $b$ (less than $v_{n}(t)$),
\begin{multline*}
\sum_{k=1}^{v_{n}(t)} \big ( {\mathbb{E}}_{k} ( X_{k} X_{k+\ell}) -
{\mathbb{E}}_{0} (X_{k}X_{k+\ell}) \big ) = \sum_{k=1}^{b} \big ( {\mathbb{E}%
}_{k} ( X_{k}X_{k+\ell} ) - {\mathbb{E}}_{0}( X_{k}X_{k+\ell} ) \big )\\
+ \sum_{k=b+1}^{v_{n}(t)} \big ( {\mathbb{E}}_{k} ( X_{k}X_{k+\ell} ) -
{\mathbb{E}}_{0} ( X_{k}X_{k+\ell} ) \big ) \, .
\end{multline*}
By the Cauchy--Schwarz inequality,
\[
\sum_{k=1}^{b} \Vert{\mathbb{E}}_{k} ( X_{k}X_{k+\ell} ) - {\mathbb{E}}_{0} (
X_{k}X_{k+\ell} ) \Vert_{1} \leq2 n^{-1} \sum_{k=1}^{b} \Vert X_{k} \Vert_{2}
\Vert X_{k+\ell} \Vert_{2} \leq2 \sum_{k=1}^{b + \ell} \Vert X_{k} \Vert
^{2}_{2} \, .
\]
Hence, for any $\varepsilon>0$,
\[
\sum_{k=1}^{b} \Vert{\mathbb{E}}_{k} ( X_{k}X_{k+\ell} ) - {\mathbb{E}}_{0} (
X_{k}X_{k+\ell} ) \Vert_{1} \leq2 \big \{ \varepsilon^{2} (b+ \ell) +
\sum_{k=1}^{b + \ell} {\mathbb{E}} \big ( X^{2}_{k} \mathbf{1}_{|X_{k}| >
\varepsilon} \big ) \big \} \, ,
\]
which converges to zero as $n \rightarrow\infty$ followed by $\varepsilon
\rightarrow0$, by taking into account condition \eqref{condLNS_g}. Now
\begin{multline*}
\sum_{k=b+1}^{v_{n}(t)} \big ( {\mathbb{E}}_{k} ( X_{k}X_{k+\ell} ) -
{\mathbb{E}}_{0} ( X_{k}X_{k+\ell} ) \big ) = \sum_{k=b+1}^{v_{n}(t)}
\big ( {\mathbb{E}}_{k} ( X_{k}X_{k+\ell}) - {\mathbb{E}}_{k-b} (
X_{k}X_{k+\ell} ) \big )\\
+ \sum_{k=b+1}^{v_{n}(t)} \big ( {\mathbb{E}}_{k-b} ( X_{k}X_{k+\ell}) -
{\mathbb{E}}_{0}( X_{k}X_{k+\ell} ) \big ) \, .
\end{multline*}
Taking into account condition \eqref{convL1carreNSrho}, we have
\[
\lim_{b \rightarrow\infty} \limsup_{n \rightarrow\infty} \Big \Vert \sum
_{k=b+1}^{v_{n}(t)} \big ( {\mathbb{E}}_{k-b} ( X_{k}X_{k+\ell} ) -
{\mathbb{E}}_{0} ( X_{k}X_{k+\ell} ) \big ) \Big \Vert_{1} = 0 \, .
\]
We show now that
\begin{equation}
\label{cequonveutpourlesP0Flo}\lim_{b \rightarrow\infty} \limsup_{n
\rightarrow\infty} \Big \Vert \sum_{k=b+1}^{v_{n}(t)} \big ( {\mathbb{E}}_{k}
( X_{k}X_{k+\ell} ) - {\mathbb{E}}_{k-b} ( X_{k}X_{k+\ell})
\big ) \Big \Vert_{1} = 0 \, .
\end{equation}
Together with the convergences proved above, this will show that
\eqref{cequonveutpourlesP0butFlo} is satisfied.

To prove \eqref{cequonveutpourlesP0Flo}, we fix a positive real $\varepsilon$
and write
\begin{multline*}
\Big \Vert \sum_{k=b+1}^{v_{n}(t)} \big ( {\mathbb{E}}_{k} ( X_{k}X_{k+\ell} )
- {\mathbb{E}}_{k-b} ( X_{k}X_{k+\ell} ) \big ) \Big \Vert_{1}\\
\leq\Big \Vert \sum_{k=b+1}^{v_{n}(t)} \big ( {\mathbb{E}}_{k} ( Y_{k,\ell
}^{\prime}) - {\mathbb{E}}_{k-b}( Y_{k,\ell}^{\prime}) \big ) \Big \Vert_{1} +
\Big \Vert \sum_{k=b+1}^{v_{n}(t)} \big ( {\mathbb{E}}_{k} ( Y_{k,\ell
}^{\prime\prime}) - {\mathbb{E}}_{k-b} ( Y_{k,\ell}^{\prime\prime})
\big ) \Big \Vert_{1}%
\end{multline*}
where
\[
Y_{k,\ell}^{\prime}= X_{k}X_{k+\ell}\mathbf{1}_{| X_{k}X_{k+\ell}|
\leq\varepsilon^{2} } \ \text{ and } \ Y_{k,\ell}^{\prime\prime}=
X_{k}X_{k+\ell}\mathbf{1}_{| X_{k}X_{k+\ell}| > \varepsilon^{2} } \, .
\]
Note now that the following inequalities are valid: for any reals $a$ and $b$
and any positive real $M$,
\begin{align}
\label{trivial20-07}|ab| \mathbf{1}_{\{|ab| >M\} } \leq2^{-1 } \big ( |a^{2}%
+b^{2}| \mathbf{1}_{\{|a^{2}+b^{2}| >2M\} } \Big ) \leq a^{2} \mathbf{1}_{\{
a^{2} >M \} } + b^{2} \mathbf{1}_{\{ b^{2} >M \} } \, .
\end{align}
Hence,
\[
\Big \Vert \sum_{k=b+1}^{v_{n}(t)} \big ( {\mathbb{E}} ( Y_{k,\ell}%
^{\prime\prime}| {\mathcal{F}}_{k}) - {\mathbb{E}} ( Y_{k,\ell}^{\prime\prime
}| {\mathcal{F}}_{k-b}) \big ) \Big \Vert_{1} \leq4 \sum_{ k=1}^{ n }
{\mathbb{E}} (X_{k}^{2} \mathbf{1}_{ |X_{k}| > \varepsilon}) \, ,
\]
which together with condition \eqref{condLNS_g} imply that
\begin{equation}
\label{cequonveutpourlesP0p1Flo}\lim_{n \rightarrow\infty} \Big \Vert \sum
_{k=b+1}^{v_{n}(t)} \big ( {\mathbb{E}}_{k} ( Y_{k,\ell}^{\prime\prime}) -
{\mathbb{E}}_{k-b}( Y_{k,\ell}^{\prime\prime}) \big ) \Big \Vert_{1} =0 \, .
\end{equation}
On another hand,
\[
\sum_{k=b+1}^{v_{n}(t)} \big ( {\mathbb{E}}_{k} ( Y_{k,\ell}^{\prime}) -
{\mathbb{E}}_{k-b} ( Y_{k,\ell}^{\prime}) \big ) = \sum_{k=b+1}^{v_{n}(t)}
\sum_{j=0}^{b-1} \mathbf{P}_{k-j} (Y_{k,\ell}^{\prime}) = \sum_{j=0}^{b-1}
\sum_{k=b+1}^{v_{n}(t)} \mathbf{P}_{k-j} ( Y_{k,\ell}^{\prime}) \, ,
\]
where we recall $\mathbf{P}_{j} ( \cdot)= {\mathbb{E}} ( \cdot| {\mathcal{F}%
}_{j}) - {\mathbb{E}} ( \cdot| {\mathcal{F}}_{j-1}) $. Since $( \mathbf{P}%
_{k-j} ( Y_{k,\ell}^{\prime}) )_{k\geq1}$ is a sequence of martingale
differences,
\[
\Big \Vert\sum_{k=b+1}^{v_{n}(t)} \big ( {\mathbb{E}}_{k} ( Y_{k,\ell}%
^{\prime}) - {\mathbb{E}}_{k-b} ( Y_{k,\ell}^{\prime}) \big ) \Big \Vert _{1 }
\leq\sum_{j=0}^{b-1} \Big \Vert \sum_{k=b+1}^{v_{n}(t)} \mathbf{P}_{k-j} (
Y_{k,\ell}^{\prime}) \Big \Vert_{2} \leq\sum_{j=0}^{b-1} \Big ( \sum
_{k=b+1}^{v_{n}(t)}\Vert\mathbf{P}_{k-j} ( Y_{k,\ell}^{\prime} ) \Vert^{2}_{2}
\Big )^{1/{2}} \, .
\]
By the Cauchy--Schwarz inequality,
\[
\Vert\mathbf{P}_{k-j} ( Y_{k,\ell}^{\prime} )\Vert^{2}_{2} \leq\Vert
{\mathbb{E}}_{k-j} ( Y_{k,\ell}^{\prime} )\Vert^{2}_{2} \leq\varepsilon^{2}
\Vert X_{k}X_{k+\ell} \Vert_{1} \leq2^{-1}\varepsilon^{2} ( \Vert X_{k}%
\Vert^{2}_{2} + \Vert X_{k+\ell} \Vert^{2}_{2} ) \, .
\]
Therefore
\[
\Big \Vert\sum_{k=b+1}^{v_{n}(t)} \big ( {\mathbb{E}}_{k} ( Y_{k,\ell}%
^{\prime}) - {\mathbb{E}}_{k-b} ( Y_{k,\ell}^{\prime}) \big ) \Big \Vert _{1 }
\leq b \, \varepsilon\sup_{n \geq1} \Big ( \sum_{k=b+1}^{n }\Vert X_{k}
\Vert^{2} _{2} \Big )^{1/{2}} \, ,
\]
which converges to zero by taking into account condition \eqref{condLNS_g} and
by letting ${\varepsilon}$ going to $0$. This last convergence together with
\eqref{cequonveutpourlesP0p1Flo} entail \eqref{cequonveutpourlesP0Flo} and
then \eqref{cequonveutpourlesP0butFlo}. This ends the proof of the
proposition. $\diamond$

\subsection{Proof of Theorem \ref{ThrhoNSalter}}

We apply Theorem \ref{fclt_MPU} to the triangular array $\{ \sigma_{n}^{-1}
X_{k,n}, 1 \leq k \leq n \}_{n \geq1}$ and the $\sigma$-algebras
${\mathcal{F}}_{k,n} = \sigma( X_{i,n}, 1 \leq i \leq k)$ for $k \geq1$ and
${\mathcal{F}}_{k,n}= \{ \emptyset, \Omega\}$ for $k \leq0$. For convenience,
we can set $X_{k,n} =0$ for $k >n$. Again, to soothe the notations, we will
omit the index $n$ involved in the variables and in the $\sigma$-algebras.

As a matter of fact, we shall first prove that under the conditions of Theorem
\ref{ThrhoNSalter}, the following reinforced version of condition
\eqref{MWcondNS2_g} is satisfied:
\begin{equation}
\lim_{m\rightarrow\infty}m^{-1/2}A(m)=0\text{ and }\lim_{m\rightarrow\infty
}\sum_{\ \ell\geq[\log_{2}(m)]}B(2^{\ell},m)=0\,.
\label{MWcondNS2_greinforced}%
\end{equation}
In order to check the conditions below, we shall apply the following
inequality, derived in Theorem 1.1 in Utev (1991). More exactly, under
\eqref{condonrho} there exists a finite positive constant $\kappa$ such that
for any positive integers $a < b$,
\begin{align}
\label{utev89_g}\Vert S_{b} - S_{a} \Vert_{2}^{2} \leq\kappa\sum_{i=a+1}^{b}
\Vert X_{i} \Vert_{2}^{2} \, .
\end{align}

The first characteristic $A^{2}(m)$ defined by \eqref{defalu_g} is then
estimated as follows. Write first the following decomposition:
\begin{equation}
\label{b1carreYkm_g}\sum_{k=0}^{n-1} \Vert{\mathbb{E}}_{k} (S_{ k +m } -
S_{k}) \Vert_{2}^{2} \, \leq2\sum_{k=0}^{n-1} \Vert{\mathbb{E}}_{k} (S_{ k +m
} - S_{k+[\sqrt{m}]}) \Vert_{2}^{2} \, + 2\sum_{k=0}^{n-1} \Vert{\mathbb{E}%
}_{k} (S_{k+[\sqrt{m}]}-S_{k}) \Vert_{2}^{2}.
\end{equation}
Note now that for any integer $k$ and any positive integers $a,b$ with $a < b
$,
\begin{multline*}
\Vert{\mathbb{E}}_{k} ( S_{k +b} - S_{k + a } ) \Vert_{2}^{2} = \mathrm{cov}
\big ( {\mathbb{E}}_{k} ( S_{k +b} - S_{k +a } ) , S_{k +b} - S_{k +a}
\big )\\
\leq\rho( a ) \Vert{\mathbb{E}}_{k} ( S_{k +b} - S_{k +a } ) \Vert_{2} \Vert
S_{k +b} - S_{k + a } \Vert_{2} \, .
\end{multline*}
Hence
\[
\Vert{\mathbb{E}}_{k} ( S_{k +b} - S_{k +a } ) \Vert_{2} \leq\rho(a ) \Vert
S_{k +b} - S_{k +a} \Vert_{2} \, ,
\]
which combined with \eqref{utev89_g} implies, under \eqref{condonrho}, that
there exists a finite positive constant $\kappa$ such that that
\begin{align}
\label{inecovrho_g}\Vert{\mathbb{E}}_{k} ( S_{k +b} - S_{k +a } ) \Vert
^{2}_{2} \leq\kappa\rho^{2} (a ) \sum_{i=k+a+1}^{k+b} \Vert X_{i} \Vert
_{2}^{2}\, .
\end{align}
Therefore, starting from \eqref{b1carreYkm_g} and taking into account
\eqref{utev89_g} and \eqref{inecovrho_g}, we get, under \eqref{condonrho} and
\eqref{condboundvar}, that
\begin{align*}
\sigma_{n}^{-2} \sum_{k=0}^{n-1} \Vert{\mathbb{E}}_{k} (S_{ k +m } - S_{k})
\Vert_{2}^{2}  &  \leq2 \kappa\sigma_{n}^{-2} \sum_{k=0}^{n-1} \rho^{2} (
[\sqrt{m}] ) \sum_{i=k +[\sqrt{m}]+1}^{k+m} \Vert X_{i} \Vert_{2}^{2} + 2
\kappa\sigma_{n}^{-2} \sum_{k=0}^{n-1} \sum_{i=k+1}^{k+[\sqrt{m}]} \Vert X_{i}
\Vert_{2}^{2}\\
&  \leq2 \kappa C \big \{ m \rho^{2} ( [\sqrt{m}] ) + \sqrt{m} \big \} \, .
\end{align*}
Hence
\begin{align}
\label{normNSrhop1}m^{-1} A^{2}(m ) \leq2 \kappa C \big \{ \rho^{2} (
[\sqrt{m}] ) + m^{-1/2} \big \} \, ,
\end{align}
which tends to zero as $m \rightarrow\infty$. This proves the first part of
assumption \eqref{MWcondNS2_greinforced}.

Next, observe that by \eqref{inecovrho_g}, under \eqref{condonrho},
\begin{align*}
\Big \Vert \frac{1}{m} \sum_{u=0}^{m-1} {\mathbb{E}} ( S_{k2^{r} + u } - S_{
(k-1) 2^{r} +u} | {\mathcal{F}}_{(k-2) 2^{r}+1} ) \Big \Vert^{2}_{2}  &
\leq\frac{1}{m} \sum_{u=0}^{m-1} \Vert{\mathbb{E}} ( S_{k2^{r} + u } - S_{
(k-1) 2^{r} +u} | {\mathcal{F}}_{(k-2) 2^{r}+1} ) \Vert^{2}_{2}\\
&  \leq\frac{\kappa}{m} \sum_{u=0}^{m-1} \rho^{2} ( 2^{r} +u-1) \sum_{i=(k-1)
2^{r} +u+1}^{k2^{r} + u } \Vert X_{i} \Vert^{2}_{2} \, .
\end{align*}
Thus, by taking into account \eqref{condboundvar}, we derive that
\begin{align*}
B^{2} ( 2^{r},m)  &  = \sup_{n \geq1} \sigma_{n}^{-2} \sum_{k=1}^{[n/2^{r}]
+1}\Big \Vert \frac{1}{m} \sum_{u=0}^{m-1} {\mathbb{E}} ( S_{k2^{r} + u } -
S_{ (k-1) 2^{r} +u} | {\mathcal{F}}_{(k-2) 2^{r}+1} ) \Big \Vert^{2}_{2}\\
&  \leq C\frac{\kappa}{m} \sum_{u=0}^{m-1} \rho^{2} ( 2^{r} +u-1) \leq C
\kappa\rho^{2} ( 2^{r} -1) \, ,
\end{align*}
where the last inequality comes from the fact that $\rho$ is non-increasing.
Taking into account \eqref{condonrho}, this shows that the second part of
\eqref{MWcondNS2_greinforced} is satisfied.

Now, we apply Proposition \ref{toverifycondSqNS_g} to verify the last
condition \eqref{AssGr3_1_g}. To do it we need to verify its assumptions
\eqref{convL1carreNSrho} and \eqref{condSqNSavecE0rho} by recalling that since
${\mathcal{F}}_{0,n}= \{ \emptyset, \Omega\}$, ${\mathbb{E}}_{0}%
(\cdot)={\mathbb{E}}(\cdot)$.

First, we notice that by the definition of the $\rho$-mixing coefficients and
the condition \eqref{condboundvar}, for any non-negative integer $\ell$,
\begin{align*}
&  \sigma_{n}^{-2}\sum_{k=b+1}^{n} \Vert{\mathbb{E}}_{k-b} ( X_{k}X_{k+\ell} )
- {\mathbb{E}} ( X_{k}X_{k+\ell} ) \Vert_{1} \leq\rho(b) \sigma_{n}^{-2}
\sum_{k=b+1}^{n} \Vert X_{k}X_{k+\ell} - {\mathbb{E}} ( X_{k}X_{k+\ell} )
\Vert_{2}\\
&  \leq\rho(b) \sigma_{n}^{-2} \sum_{k=b+1}^{n} \Vert X_{k} \Vert_{2} \Vert
X_{k+\ell} \Vert_{2} \leq\rho(b) \sigma_{n}^{-2} \Big (\sum_{k=1}^{n+ \ell}
\Vert X_{k} \Vert_{2}^{2} \Big ) \leq\rho(b) C\to0 \mbox{ as } b\to\infty,
\end{align*}
which proves the first assumption \eqref{convL1carreNSrho}.

\medskip

To end the proof of the theorem, it remains to prove that
\eqref{condSqNSavecE0rho} holds. Note that since we have proved that condition
\eqref{MWcondNS2_greinforced} is satisfied, a careful analysis of the proof of
Lemma \ref{Gen_fclt} reveals that, setting $D_{\ell}^{m}=m^{-1}\sum_{i={0}%
}^{m-1}{\mathbf{P}}_{\ell}(S_{{\ell}+i})$ and $\theta_{{\ell}}^{m}=m^{-1}%
\sum_{i=1}^{m-1}{\mathbb{E}}_{{\ell}}(X_{{\ell+1}}+\ldots+X_{{\ell}+i})$,
\begin{equation}
\lim_{m\rightarrow\infty}\limsup_{n\rightarrow\infty}\sigma_{n}^{-2}  \sup_{t\in\lbrack0,1]}
\Big |\sum_{\ell=1}^{v_{n}(t)}\Big ({\mathbb{E}}(X_{{\ell}}^{2}+2X_{\ell
}\theta_{{\ell}}^{m})-{\mathbb{E}}(D_{\ell}^{m})^{2}%
\Big )\Big |=0\,,\label{condSqNSavecE0rhoeasier}%
\end{equation}
and
\begin{equation}
\lim_{m \rightarrow \infty}\limsup_{n\rightarrow\infty}\sigma_{n}^{-2} \Big \Vert  \sup_{t\in\lbrack0,1]} \Big |\sum_{i=1}^{v_{n}(t)}X_{i,n}%
-\sum_{\ell=1}^{v_{n}(t)}D_{{\ell}}^{m} \Big |\Big \Vert_{2}^2=0\,.\label{mart approx(2)avecm}%
\end{equation}
Taking into account \eqref{condSqNSavecE0rhoeasier}, to prove that
\eqref{condSqNSavecE0rho} holds, we then need to show that, for any
$t\in\lbrack0,1]$,
\begin{equation}
\lim_{m\rightarrow\infty}\limsup_{n\rightarrow\infty}\big |\sigma_{n}^{-2}%
\sum_{\ell=1}^{v_{n}(t)}{\mathbb{E}}(D_{\ell}^{m})^{2}%
-t\big |=0\,.\label{condSqNSavecE0rhointD}%
\end{equation}
But, since $\sum_{\ell=1}^{v_{n}(t)}{\mathbb{E}}(D_{\ell}^{m})^{2}=\Vert
\sum_{\ell=1}^{v_{n}(t)}D_{\ell}^{m}\Vert_{2}^{2}$, by
taking into account \eqref{mart approx(2)avecm}, the convergence
\eqref{condSqNSavecE0rhointD} will follow if one can prove that, for any
$t\in\lbrack0,1]$,
\begin{equation}
\sigma_{n}^{-2}{\mathbb{E}}\big (S_{v_{n}(t)}^{2}\big )\rightarrow
t\,,\,\text{ as $n\rightarrow\infty$}\,.\label{constatationtrivial_g}%
\end{equation}
With this aim, we note that since
\[
\Big \Vert\sum_{k=1}^{v_{n}(t)}X_{k,n}\Big \Vert_{2}\leq\Big \Vert\sum
_{k=1}^{v_{n}(t)-1}X_{k,n}\Big \Vert_{2}+\Vert X_{v_{n}(t),n}\Vert_{2}\,,
\]
by definition of $v_{n}(t)$, we have
\[
\sqrt{t}\leq\frac{\Vert\sum_{k=1}^{v_{n}(t)}X_{k,n}\Vert_{2}}{\sigma_{n}}%
\leq\sqrt{t}+\frac{\Vert X_{v_{n}(t),n}\Vert_{2}}{\sigma_{n}}\,.
\]
This implies \eqref{constatationtrivial_g} by noticing that the Lindeberg
condition \eqref{condLNSavecsigma} implies that $\lim_{n\rightarrow\infty
}\frac{\Vert X_{v_{n}(t),n}\Vert_{2}}{\sigma_{n}}=0$. The proof of Theorem
\ref{ThrhoNSalter} is complete. $\diamond$

\subsection{ Proof of Corollary \ref{ThrhoNS}}

By taking $v_{n}(t)=[nt]$, we need to ensure that
\eqref{constatationtrivial_g} holds, which is straightforward since we assume
that $\sigma_{n}^{2}=n h(n)$ where $h$ is a slowly varying function at
infinity. $\diamond$

\subsection{Proof of Corollary \ref{corforlinear}}

We note first that, for any $k \geq1$,
\[
\Vert X_{k} \Vert_{2} \leq2 \Vert f(Y_{k}) - f(0) \Vert_{2} \leq2 \Vert c(
|Y_{k} | ) \Vert_{2} \leq2 c (\Vert Y_{k} \Vert_{2}) \, ,
\]
where, since $c$ is non-decreasing and concave, we have used Lemma 5.1 in
Dedecker (2008). Therefore by \eqref{Kac}
\[
\sup_{k \geq1} \Vert X_{k} \Vert_{2} \leq2 c \Big ( \sigma_{\varepsilon}
\sum_{i \geq0} |a_{i}| \Big ) < \infty\, .
\]
This proves the first part of \eqref{condLNS}. Now, to prove the second part of \eqref{condLNS}, it
suffices to show that, for any $\varepsilon>0$,
\begin{equation}
\lim_{n\rightarrow\infty}\frac{1}{n}\sum_{k=1}^{n} {\mathbb{E}} \big ( |
f(Y_{k}) - f(0) |^{2} \mathbf{1}_{\{ |f(Y_{k}) - f(0)|>\varepsilon\sqrt{n})
\}} \big ) =0 \, . \label{condLNSlinearproc}%
\end{equation}
With this aim, we set $C= \sum_{i \geq0} |a_{i}| $ and let $K$ be a positive
integer. We denote by $Y^{\prime}_{k}= \sum_{i \geq0} a_{i} \varepsilon
^{\prime}_{k-i}$ where $\varepsilon_{i}^{\prime}= \varepsilon_{i}
\mathbf{1}_{|\varepsilon_{i} | \leq K}$. Let also $\varepsilon_{i}%
^{\prime\prime}= \varepsilon_{i} \mathbf{1}_{|\varepsilon_{i} |> K}$. Using
the fact that for positive reals $a, b$ and $\varepsilon$, $(a+b)^{2}
\mathbf{1}_{\{ a+b >2 \varepsilon\}} \leq4 a^{2} \mathbf{1}_{\{ a >
\varepsilon\}} + 4b^{2} \mathbf{1}_{\{ b > \varepsilon\}}$, we get that, for
any $\varepsilon>0$,
\begin{multline*}
{\mathbb{E}} \big ( | f(Y_{k}) - f(0) |^{2} \mathbf{1}_{\{ |f(Y_{k}) -
f(0)|>2\varepsilon\sqrt{n}) \}} \big ) \leq4 {\mathbb{E}} \big ( | f(Y_{k}) -
f(Y_{k}^{\prime }) |^2\mathbf{1}_{\{ |f(Y_{k}) - f(Y_{k}^{\prime})|>\varepsilon
\sqrt{n}) \}} \big )\\
+ 4 {\mathbb{E}} \big ( | f(Y^{\prime}_{k}) - f(0) |^{2} \mathbf{1}_{\{
|f(Y^{\prime}_{k}) - f(0)|>\varepsilon\sqrt{n}) \}} \big )\\
\leq4 \Big \Vert c \Big ( \sum_{i \geq0} | a_{i} \varepsilon^{\prime\prime
}_{k-i}| \Big ) \Big \Vert_{2}^{2} + 4 c^{2} \big ( K C \big ) \mathbf{1}_{\{
c (KC)>\varepsilon\sqrt{n}) \}} \, .
\end{multline*}
The last term in the right-hand side converges to zero as $n \rightarrow
\infty$. Next, since $c$ is non-decreasing and concave, Lemma 5.1 in Dedecker
(2008) gives
\[
\frac{1}{n}\sum_{k=1}^{n} \Big \Vert c \Big ( \sum_{i \geq0} | a_{i}
\varepsilon^{\prime\prime}_{k-i}| \Big ) \Big \Vert_{2}^{2} \leq\frac{1}%
{n}\sum_{k=1}^{n} c^{2} \Big ( \sum_{i \geq0} | a_{i} | \Vert\varepsilon
^{\prime\prime}_{k-i} \Vert_{2} \Big ) \leq c^{2} \Big ( \sup_{k
\in{\mathbb{Z}}} \Vert\varepsilon^{\prime\prime}_{k} \Vert_{2} \sum_{i \geq0}
| a_{i} | \Big ) \, .
\]
But, since $(\varepsilon^{2}_{i})_{i \in{{\mathbb{Z}}}}$ is an uniformly
integrable family, $\limsup_{K \rightarrow\infty} \sup_{k \in{\mathbb{Z}}}
\Vert\varepsilon^{\prime\prime}_{k} \Vert_{2} =0 $. Together with the fact
that $\lim_{x \rightarrow0} c(x) =0$, this proves that
\[
\lim_{K\rightarrow\infty}\limsup_{n\rightarrow\infty} \frac{1}{n}\sum
_{k=1}^{n} \Big \Vert c \Big ( \sum_{i \geq0} | a_{i} \varepsilon
^{\prime\prime}_{k-i}| \Big ) \Big \Vert_{2}^{2} =0 \, ,
\]
ending the proof of \eqref{condLNSlinearproc} and then of \eqref{condLNS}.

Let us consider now the following choice of $({\mathcal{F}}_{i})_{i \geq0}$:
${\mathcal{F}}_{0} = \{ \emptyset, \Omega\}$ and ${\mathcal{F}}_{i}= \sigma(
X_{1}, \ldots, X_{i} )$, for $i \geq1$. If one can prove that conditions
\eqref{MWcondNScasereinforced} and \eqref{convL1carreNSrho} are satisfied and
also that, for any $t \in[0,1]$,
\begin{equation}
\lim_{m\rightarrow\infty}\limsup_{n\rightarrow\infty} \frac{1}{\sigma_{n}^{2}}
\Big | \sum_{k=1}^{[nt]} \Big \{ {\mathbb{E}}(X_{k}^{2})+ 2 {\mathbb{E}}
(X_{k}\theta^{m}_{k} ) \Big \} - t \Big | =0\, , \label{condSqNSavecE0rhobis}%
\end{equation}
then the corollary will follow by applying Theorem
\ref{CorollaryCLTNSforappli} and by taking into account Proposition
\ref{toverifycondSqNS_g}.

\smallskip

To prove that \eqref{MWcondNScasereinforced} holds, we set ${\mathbb{E}%
}_{\varepsilon}$ the expectation with respect to $\varepsilon:= (\varepsilon
_{i})_{i \in{{\mathbb{Z}}}}$ and note that since ${\mathcal{F}}_{i}
\subset{\mathcal{F}}_{\varepsilon, i}$ where ${\mathcal{F}}_{\varepsilon, i} =
\sigma( \varepsilon_{k}, k \leq i )$, for any $i \geq0$,
\begin{align}
\label{boundorder1}\Vert{\mathbb{E}}(X_{ k+i } |{\mathcal{F}}_{i})\Vert_{2}
\leq\Vert{\mathbb{E}}(X_{ k+i } |{\mathcal{F}}_{\varepsilon, i} )\Vert_{2} \,
.
\end{align}
For any $i \geq0$,
\begin{multline*}
\big | {\mathbb{E}}(X_{ k+i } |{\mathcal{F}}_{\varepsilon, i} ) \big | =
\Big | {\mathbb{E}}_{\varepsilon} \Big ( f \Big ( \sum_{\ell=0}^{k-1} a_{\ell
}\varepsilon^{\prime}_{k +i - \ell} + \sum_{\ell\geq k} a_{\ell}\varepsilon_{k
+i - \ell} \Big )\Big ) - {\mathbb{E}}_{\varepsilon} \Big ( f \Big ( \sum
_{\ell=0}^{k-1} a_{\ell}\varepsilon^{\prime}_{k +i - \ell} + \sum_{\ell\geq k}
a_{\ell}\varepsilon^{\prime}_{k +i - \ell} \Big )\Big ) \Big | \, ,
\end{multline*}
where $(\varepsilon^{\prime}_{i})_{i \in{{\mathbb{Z}}}}$ is an independent
copy of $(\varepsilon_{i})_{i \in{{\mathbb{Z}}}}$. Therefore, using Lemma 5.1
in Dedecker (2008),
\[
\Vert{\mathbb{E}}(X_{ k+i } |{\mathcal{F}}_{i})\Vert_{2} \leq\Big \Vert c
\Big ( \sum_{\ell\geq k} |a_{\ell}| | \varepsilon_{k +i - \ell} -
\varepsilon^{\prime}_{k +i - \ell} | \Big ) \Big \Vert_{2} \leq c \Big ( 2
\sigma_{\varepsilon} \sum_{\ell\geq k} |a_{\ell}| \Big ) \, ,
\]
proving that \eqref{MWcondNScasereinforced} holds under \eqref{Kac}. We prove
now that \eqref{convL1carreNSrho} is satisfied. With this aim we recall that
${\mathcal{F}}_{0}$ is the trivial $\sigma$-field, and we first write that for
any non-negative integer $k,j$, and $n$,
\begin{multline}
\label{decproduct1}{\mathbb{E}}_{k}(X_{k+n} X_{j +n} ) - {\mathbb{E}}(X_{k+n}
X_{j +n} ) = {\mathbb{E}}_{k} (f(Y_{k+n}) f(Y_{j +n} )) - {\mathbb{E}}
(f(Y_{k+n}) f(Y_{j +n} ))\\
- {\mathbb{E}}(f(Y_{k+n}) ) {\mathbb{E}}_{k} ( f(Y_{j +n} ) - {\mathbb{E}}(
f(Y_{j +n} )) ) - {\mathbb{E}}(f(Y_{j+n}) ) {\mathbb{E}}_{k} ( f(Y_{k+n} ) -
{\mathbb{E}}( f(Y_{k +n} )) ) \, .
\end{multline}
Now, for $j \geq k$,
\begin{multline*}
{\mathbb{E}}(f(Y_{k+n}) f(Y_{j +n} ) | {\mathcal{F}}_{\varepsilon, k} ) -
{\mathbb{E}} (f(Y_{k+n}) f(Y_{j +n} ))\\
= {\mathbb{E}}_{\varepsilon} \Big ( f \Big ( \sum_{\ell=0}^{n-1} a_{\ell
}\varepsilon^{\prime}_{k +n - \ell} + \sum_{\ell\geq n} a_{\ell}\varepsilon_{k
+n - \ell} \Big ) f \Big ( \sum_{\ell=0}^{n-1 +j-k } a_{\ell}\varepsilon
^{\prime}_{j +n - \ell} + \sum_{\ell\geq n +j-k} a_{\ell}\varepsilon_{j +n -
\ell} \Big ) \Big )\\
- {\mathbb{E}}_{\varepsilon} \Big ( f \Big ( \sum_{\ell=0}^{n-1} a_{\ell
}\varepsilon^{\prime}_{k +n - \ell} + \sum_{\ell\geq n} a_{\ell}%
\varepsilon^{\prime}_{k +n - \ell} \Big ) f \Big ( \sum_{\ell=0}^{n-1 +j-k }
a_{\ell}\varepsilon^{\prime}_{j +n - \ell} + \sum_{\ell\geq n +j-k} a_{\ell
}\varepsilon^{\prime}_{j +n - \ell} \Big ) \Big ) \, .
\end{multline*}
Setting $Z_{k} = \sum_{\ell=0}^{n-1} a_{\ell}\varepsilon^{\prime}_{k +n -
\ell} + \sum_{\ell\geq n} a_{\ell}\varepsilon_{k +n - \ell} $, $Z^{\prime}_{k}
= \sum_{\ell\geq0} a_{\ell}\varepsilon^{\prime}_{k +n - \ell} $ and $U_{k,j} =
\sum_{\ell=0}^{n-1 +j-k } a_{\ell}\varepsilon^{\prime}_{j +n - \ell} +
\sum_{\ell\geq n +j-k} a_{\ell}\varepsilon_{j +n - \ell} $, we get
\begin{multline*}
\Vert{\mathbb{E}}(f(Y_{k+n}) f(Y_{j +n} ) | {\mathcal{F}}_{\varepsilon, k} ) -
{\mathbb{E}} (f(Y_{k+n}) f(Y_{j +n} )) \Vert_{1}
\leq\Vert f ( Z_{k} ) f ( U_{k,j}) - f ( Z^{\prime}_{k} ) f ( Z^{\prime}_{j})
\Vert_{1}\\
\leq\Vert f ( Z_{k} ) - f ( Z^{\prime}_{k} ) \Vert_{2} \Vert f ( U_{k,j})
\Vert_{2} + \Vert f ( U_{k,j}) - f ( Z^{\prime}_{j}) \Vert_{2} \Vert f (
Z^{\prime}_{k} ) \Vert_{2} \, .
\end{multline*}
By using Lemma 5.1 in Dedecker (2008),
\[
\Vert f ( Z_{k} ) - f ( Z^{\prime}_{k} ) \Vert_{2} \leq\Big \Vert c
\Big ( \sum_{\ell\geq n} |a_{\ell}| | \varepsilon_{k +n - \ell} -
\varepsilon^{\prime}_{k +n - \ell} | \Big ) \Big \Vert_{2} \leq c \Big ( 2
\sigma_{\varepsilon} \sum_{\ell\geq n} |a_{\ell}| \Big ) \, .
\]
Similarly
\[
\Vert f ( U_{k,j}) - f ( Z^{\prime}_{j}) \Vert_{2} \leq c \Big ( 2
\sigma_{\varepsilon} \sum_{\ell\geq n +j-k } |a_{\ell}| \Big ) \, .
\]
Moreover, using Lemma 5.1 in Dedecker (2008) again we have
\[
\Vert f ( U_{k,j}) \Vert_{2} \leq f(0) + c (\Vert U_{k,j} \Vert_{2}) \, .
\]
Therefore, by \eqref{Kac},
\[
\sup_{j \geq k \geq0} \Vert f ( U_{k,j}) \Vert_{2} \leq f(0) + c
\Big ( \sigma_{\varepsilon} \sum_{i \geq0} |a_{i} | \Big ) < \infty\, .
\]
Similarly
\[
\sup_{k \geq0} \Vert f ( Z^{\prime}_{k} ) \Vert_{2} \leq f(0) + c
\Big ( \sigma_{\varepsilon} \sum_{i \geq0} |a_{i} | \Big ) < \infty\, .
\]
So, overall, and since $\lim_{x \rightarrow0} c(x) =0$, we get
\begin{align}
\label{decproduct2}\lim_{n \rightarrow\infty} \sup_{ j \geq k \geq0 }
\Vert{\mathbb{E}}_{k}(f(Y_{k+n}) f(Y_{j +n} ) - {\mathbb{E}} (f(Y_{k+n})
f(Y_{j +n} )) \Vert_{1} =0 \, .
\end{align}
On another hand, by using \eqref{boundorder1}, for any non-negative integers
$k,j,n$,
\[
\Vert{\mathbb{E}}(f(Y_{k+n}) ) {\mathbb{E}}_{k} ( f(Y_{j +n} ) - {\mathbb{E}}(
f(Y_{j +n} )) ) \Vert_{1} \leq\sup_{\ell\geq0} \Vert f(Y_{\ell}) \Vert_{2} c
\Big ( 2 \sigma_{\varepsilon} \sum_{\ell\geq n+j - k} |a_{\ell}| \Big ) \, .
\]
Since, as before, $\sup_{\ell\geq0} \Vert f(Y_{\ell}) \Vert_{2} \leq f(0) + c
\Big (   \sigma_{\varepsilon} \sum_{i \geq0} a_{i} \Big )  < \infty$ and
$\lim_{x \rightarrow0} c(x) =0$, it follows that
\begin{align}
\label{decproduct3}\lim_{n \rightarrow\infty} \sup_{ j \geq k \geq0 }
\Vert{\mathbb{E}}(f(Y_{k+n}) ) {\mathbb{E}}_{k} ( f(Y_{j +n} ) - {\mathbb{E}}(
f(Y_{j +n} )) ) \Vert_{1} =0 \, .
\end{align}
Similarly,
\[
\Vert{\mathbb{E}}(f(Y_{j+n}) ) {\mathbb{E}}_{k} ( f(Y_{k+n} ) - {\mathbb{E}}(
f(Y_{k +n} )) ) \Vert_{1} \leq\sup_{\ell\geq0} \Vert f(Y_{\ell}) \Vert_{2} c
\Big ( 2 \sigma_{\varepsilon} \sum_{\ell\geq n} |a_{\ell}| \Big ) \, ,
\]
and therefore
\begin{align}
\label{decproduct3bis}\lim_{n \rightarrow\infty} \sup_{ j \geq k \geq0 }
\Vert{\mathbb{E}}(f(Y_{j+n}) ) {\mathbb{E}}_{k} ( f(Y_{k+n} ) - {\mathbb{E}}(
f(Y_{k +n} )) ) \Vert_{1} =0 \, .
\end{align}
Starting from \eqref{decproduct1} and taking into account
\eqref{decproduct2}-\eqref{decproduct3bis}, the convergence
\eqref{convL1carreNSrho} follows since we assumed that $\sigma_{n}^{2} = n
h(n)$ where $h(n)$ is a slowly varying function at infinity such that
$\liminf_{n \rightarrow\infty} h(n) >0$.

\medskip

We turn now to the proof of \eqref{condSqNSavecE0rhobis}. With this aim, note
first that since condition \eqref{MWcondNScasereinforced} is satisfied and
$\sigma_{n}^{2} = n h(n)$ where $h(n)$ is a slowly varying function at
infinity such that $\liminf_{n \rightarrow\infty} h(n) >0$, condition
\eqref{MWcondNS2_greinforced} holds. Now as quoted in the proof of Theorem
\ref{ThrhoNSalter}, if the Lindeberg-type condition \eqref{condLNS_g} and
condition \eqref{MWcondNS2_greinforced} are both satisfied, then to prove
\eqref{condSqNSavecE0rhobis} it is enough to show that
\eqref{constatationtrivial_g} holds (here with $v_{n}(t) =[nt]$). This comes
obviously from the fact that we assumed that $\sigma_{n}^{2} = n h(n)$ where
$h(n)$ is a slowly varying function at infinity. This ends the proof of the
corollary. $\diamond$

\subsection{Proof of Corollary \ref{corquenched}}

For any integrable random variable $f$ from $\Omega$ to ${\mathbb{R}}$ we
write $K(f)=P_{T|{\mathcal{F}}_{0}}(f)$. Since ${\mathbb{P}}$ is invariant by
$T$, for any integer $k$, a regular version $P_{T|{\mathcal{F}}_{k}}$ of $T$
given ${\mathcal{F}}_{k}$ is then obtained $via$ $P_{T|{\mathcal{F}}_{k}%
}(f)=K(f\circ T^{-k})\circ T^{k}$. With these notations, for any positive
integer $\ell$, ${\mathbb{E}}(f\circ T^{\ell}|{\mathcal{F}}_{0})=K^{\ell}(f)$.
We denote
\[
{\mathcal{M}}_{2^{r}}(|f|)=\sup_{n\geq1}\frac{1}{n}\sum_{k=0}^{n-1}K^{k2^{r}%
}(|f|)\, .
\]

Applying Corollary \ref{comment1}, Corollary \ref{corquenched} will follow if
one can prove that, with probability one,
\begin{equation} \label{cond1quenched}
\sup_{n\geq1}n^{-1}\sum_{j=1}^{n}{\mathbb{E}}_{0}(X_{j}^{2})\leq C<\infty\,,
\end{equation}
\begin{equation}
\lim_{n\rightarrow\infty}\frac{1}{n}\sum_{k=1}^{n}{\mathbb{E}}_{0}\{X_{k}%
^{2}I(|X_{k}|>\varepsilon\sqrt{n})\}=0\,,\,\text{ for any $\varepsilon>0$%
}\,,\label{condLNSquenched}%
\end{equation}
there exists a constant $c^{2}$ such that, for any $t\in\lbrack0,1]$ and any
$\varepsilon>0$,
\begin{equation}
\liminf_{m\rightarrow\infty}\limsup_{n\rightarrow\infty}{\mathbb{P}}%
_{0}\Big (\Big |\frac{1}{n}\sum_{k=1}^{[nt]}\big (X_{k}^{2}+\frac{2}{m}%
X_{k}\sum_{i=1}^{m-1}{\mathbb{E}}_{k}(S_{k+i}-S_{k})\big )-tc^{2}%
\Big |>\varepsilon\Big )=0\,,\label{condSqNSquenched}%
\end{equation}%
\begin{equation}
\sum_{\ell\geq0}2^{-\ell/2}{\mathcal{M}}_{1}^{1/2}(|{\mathbb{E}}%
_{0}(S_{2^{\ell}})|^{2})<\infty \, , \label{MWcondNS2quenched1}%
\end{equation}
and
\begin{equation}
\liminf_{j\rightarrow\infty}\sum_{\ \ell\geq j}2^{-\ell/2}{\mathcal{M}%
}_{2^{\ell}}^{1/2}\big (\big |2^{-j}\sum_{u=0}^{2^{j}-1}{\mathbb{E}}%
_{-2^{\ell}+1}(S_{2^{\ell}}\circ T^{u})\big |^{2}%
\big )=0\,.\label{MWcondNS2quenched2}%
\end{equation}

Let introduce the weak ${\mathbb{L}}^{2}$-spaces
\[
{\mathbb{L}}^{2,w}:=\{f\in{\mathbb{L}}^{1}~:~\sup_{\lambda>0}\lambda
^{2}{\mathbb{P}}\{|f|\geq\lambda\}<\infty\}\,.
\]
Recall that, when $p>1$, there exists a norm $\Vert\cdot\Vert_{2,w}$ on
${\mathbb{L}}^{2,w}$ that makes ${\mathbb{L}}^{2,w}$ a Banach space and which
is equivalent to the "pseudo"-norm $(\sup_{\lambda>0}\lambda^{2}{\mathbb{P}%
}\{|f|\geq\lambda\})^{1/2}$. Moreover, by the Dunford--Schwartz (or Hopf)
ergodic theorem (see, e.g., Krengel (1985), Lemma 6.1, page 51, and Corollary
3.8, page 131), there exists $C>0$ and such that for every $f\in{\mathbb{L}%
}^{2}$ and any non-negative integer $\ell$,
\begin{equation}
\Vert({\mathcal{M}}_{2^{\ell}}(|f|^{2}))^{1/2}\Vert_{2,w}\leq C\Vert
f\Vert_{2}\,.\label{inemaxmaxnorm}%
\end{equation}
With the help of \eqref{inemaxmaxnorm}, it is then easy to see that
\eqref{MWcondNS2quenched1} and \eqref{MWcondNS2quenched2} are satisfied under \eqref{MWcond}.

\smallskip

To prove \eqref{cond1quenched} and \eqref{condLNSquenched}, it suffices to
apply, for instance, Lemma 7.1 in Dedecker \textit{et al.} (2014).

\smallskip

It remains to prove that \eqref{condSqNSquenched} is satisfied. Since, under
\eqref{MWcond}, $\lim_{m\rightarrow\infty}m^{-1/2}{\mathbb{E}}(S_{m}%
^{2})=c^{2}$, by the ergodic theorem and the proof of Corollary \ref{rem2},
\[
\lim_{m\rightarrow\infty}\lim_{n\rightarrow\infty}\Big |\frac{1}{n}\sum
_{k=1}^{[nt]}\big (X_{k}^{2}+\frac{2}{m}X_{k}\sum_{i=1}^{m-1}{\mathbb{E}}%
_{k}(S_{k+i}-S_{k})\big )-tc^{2}\Big |=0\,,\,\text{almost surely.}%
\]
This proves \eqref{condSqNSquenched} by taking into account the properties of
the conditional expectation (see, e.g., Theorem 34.3, item (v) in Billingsley
(1995)). This ends the proof of the corollary. $\diamond$

\subsection{Proof of Corollary \ref{Gaussian}}

The proof is based on the fact that, by Theorem 27.5 in Bradley (2007), the
underlying stationary Gaussian sequence $(Y_{k})_{k\in{\mathbb{Z}}}$ is
$\rho-$mixing satisfying condition (\ref{condonrho}). By the definition of the
$\rho-$mixing coefficients we can easily derive that $(X_{k})_{k\in
{\mathbb{Z}}}$ is also a $\rho-$mixing sequence satisfying the rate condition
(\ref{condonrho}). The result follows by Corollary \ref{corquenchedro}.
$\diamond$

\subsection{Proof of Corollary \ref{appliRWRTSstat}} In all the proof, we recall that the chain starts from the origin. To soothe the notation, we will often write ${\mathbb P}$ for ${\mathbb P}_{\phi_0=0}$ when no confusion is possible.

We first note that $c^2$ is well defined. Indeed, the series $ \sum_{m \geq 1} \E ( \zeta_0 \zeta_m ) \sum_{j=1}^m (P^j)_{0,m-j}$ is absolutely convergent. To see this, we note that 
since $\pi$ is the stationary distribution, we have that for any non-negative integers $m$ and $i$, $\pi_m = \sum_{j \geq 0}{\pi}_j (P^{i})_{j,m}$. Therefore $ (P^j)_{0,m-j} \leq \pi_0^{-1} \pi_{m-j}$ which entails that
\[
 \sum_{m \geq 1} \big | \E ( \zeta_0 \zeta_m )  \big | \sum_{j=1}^m (P^j)_{0,m-j}  \leq  \sum_{m \geq 1} \big | \E ( \zeta_0 \zeta_m )  \big |   \sum_{j\geq 1} j p_j = \E ( \tau) \sum_{m \geq 1} \big | \E ( \zeta_0 \zeta_m )  \big |  \, .
\]
So the absolute convergence of the series $ \sum_{m \geq 1} \E ( \zeta_0 \zeta_m ) \sum_{j=1}^m (P^j)_{0,m-j}$ comes from the fact that $ \sum_{m \geq 1} \big | \E ( \zeta_0 \zeta_m )  \big |  < \infty$ which holds under 
the first part of condition \eqref{condARWRTSpart1stat}. Indeed, the first part of condition \eqref{condARWRTSpart1stat} implies that 
$\sum_{k \geq 0}  \Vert  \rE ( \zeta_k | {\mathcal G}_0 ) -\rE ( \zeta_k | {\mathcal G}_{-1} )  \Vert_2 < \infty$ (see Corollary 2 in \cite{PU2}) which in turn implies $ \sum_{m \geq 1} \big | \E ( \zeta_0 \zeta_m )  \big |  < \infty$. 

\medskip

The invariance principle  will be proved by an application of Theorem  \ref{CorollaryCLTNSforappli} together with Proposition \ref{toverifycondSqNS_g}. In order to apply it we need to introduce a non-decreasing filtration $({\mathcal F}_k)_{k \geq 0}$ for which the process $(X_k)_{k \geq 1}$ is adapted. Let 
\[
{\mathcal A}= \sigma ( Y_k, k \geq 0)  
\]
and for any $k \in {\mathbb Z}$, we define
\begin{equation*}  
\mathcal{F}_k = \sigma(\mathcal{A}; X_j, 1\leq j\leq k) \, .
\end{equation*}
This filtration encodes all history of the Markov chain, and all history up to time $k$ of the sampled time  scenery
$\zeta_{Y_1},\zeta_{Y_2},\dots, \zeta_{Y_k}$. 
Note that by the above definition we have: $\mathcal{F}_{j} = \mathcal{A}$, for all $j\leq 0$. 

\medskip

We first notice that by Condition ($A_1$) and  the independence between the time scenery and the Markov chain, \eqref{condLNS} is satisfied. It remains to show that conditions \eqref{MWcondNScasereinforced},   \eqref{convL1carreNSrho} and \eqref{condSqNSavecE0rho} (with $v_n(t)=[nt]$ and $X_{k,n}=X_k/{\sqrt n}$) hold. With this aim, we start by introducing some notations.

\smallskip

Let us first define a sequence of stopping times $\nu_1, \nu_2,\dots$ as the times of return to the origin,
$$\nu_1 \eqdef \inf\{k>0 : \phi_k = 0\}, \quad\nu_{n+1} \eqdef \inf\{k> \nu_{n}: \phi_k = 0\}, \quad n\geq 1 \, .$$
Define also  the random times $\tau_1, \tau_2, \dots$ as the interarrival times to the origin,
i.e. the times between successive returns to the origin, 
$$\tau_1 = \nu_1, \quad \tau_{n+1} = \nu_{n+1} - \nu_n, \quad n\geq 1 \, ,$$
and the renewal process
\bbb \label{defrenewproc} 
N_k = \max\{j \geq 0: \nu_j \leq k\} +1 \, ,
\eee
which basically counts the number of visits to the origin up to time $k$, including the renewal at $0$. Let us start with the two following facts that will be useful to prove that conditions \eqref{MWcondNScasereinforced} and  \eqref{convL1carreNSrho} are satisfied. Their proof will be given later. 

\begin{Fact} \label{computationPrequalRWRTS} Let $k$ and $m$ be in ${\mathbb N}$ and such that $m \leq k$. Then 
\[
\rP(N_{k-m} = N_k)  \leq \rP(\tau > m) \, .
\]
\end{Fact}

\begin{Fact} \label{computationofthecondexpectationstat} For any non-negative integer $k$, let 
\[
b(k) = \Vert \E ( \zeta_k | {\mathcal G}_0 ) \Vert_2 \, \text{ and } \, c(k) =   \sup_{j\geq i \geq k} \Vert \E ( \zeta_i \zeta_j  | {\mathcal G}_0 ) - \E ( \zeta_i \zeta_j )   \Vert_{1 } \, .
\]
Then
\[
\sup_{k \geq 0} \Vert \E ( X_{k+m} | {\mathcal G}_k ) \Vert^2_{2} \leq b^2 ([m/2]) + b^2(0) \rP(\tau > [m/2])
\]
and
\[
\sup_{k,\ell \geq 0} \Vert \E ( X_{k+m} X_{k+m+\ell}| {\mathcal F}_k )  -\E_{\mathcal A} ( X_{k+m} X_{k+m+\ell})  \Vert_{1} \leq c ([m/2]) + c(0) \rP(\tau > [m/2]) \, .
\]
\end{Fact}

The second inequality of Fact \ref{computationofthecondexpectationstat} together with the second  part of condition \eqref{condARWRTSpart1stat} and the fact that $\E (\tau) < \infty$ implies that condition \eqref{convL1carreNSrho} is satisfied. On another hand, by the first inequality of Fact \ref{computationofthecondexpectationstat} and taking into account the first part of condition \eqref{condARWRTSpart1stat}, it follows that  condition \eqref{MWcondNScasereinforced}  will be satisfied provided that 
$\sum_{m \geq 1} m^{-1/2} \sqrt{\rP(\tau > [m/2])} < \infty$. This last condition obviously holds since ${\mathbb E} (\tau^2) < \infty$ 
is equivalent to $\sum_{i \geq 1}  i \,  {\mathbb P}(\tau > i) < \infty $.

\smallskip

It remains to prove that condition \eqref{condSqNSavecE0rho}  holds. Note that it will be satisfied if one can prove that, for any $\varepsilon >0$ and any $t \in [0,1]$, 
\bbb \label{toprovealternativesecondcond}
\lim_{m\rightarrow \infty }\lim \sup_{n\rightarrow \infty }{\mathbb P} \Big (  \frac{1}{m}  \sum_{K=1}^{m-1} \Big |\frac{1}{n}%
\sum_{k=1}^{[nt]} \Big \{  \E_{\A} (X_{k}^{2})+ 2  \sum_{\ell=1}^K  \E_{\A} (X_k X_{k + \ell}  ) \Big \} - t c^{2} \Big |>\varepsilon \Big ) =0  \, .
\eee
Note that for any integer $\ell$ in $[0,K]$ where $K$ is a fixed integer, because of the independence between the time scenery and the Markov chain,  by stationarity of the random time scenery,
\begin{multline*}
\E_{\mathcal A} ( X_k X_{k+\ell}) = \sum_{m \geq \ell} {\bf 1}_{\phi_k = m} \E ( \zeta^2_{k+m}) + \sum_{i=1}^{\ell}  \sum_{m \geq 0} {\bf 1}_{\phi_k = \ell-i}{\bf 1}_{\phi_{k+\ell} = m} \E ( \zeta_{k+\ell-i} \zeta_{k+\ell+ m})  \\
=\E ( \zeta^2_{0})  {\bf 1}_{\phi_k \geq \ell } +  \sum_{i=1}^{\ell}  \sum_{m \geq 0} {\bf 1}_{\phi_k = \ell-i}{\bf 1}_{\phi_{k+\ell} = m} \E ( \zeta_{0} \zeta_{ m+i}) \, .
\end{multline*}
It follows that 
\begin{multline}  \label{toprovealternativesecondcondp1}
 \Big | \frac{1}{n}\sum_{k=1}^{[nt]} \big (  \E_{\mathcal A}(X_k^{2}) + 2 \sum_{\ell=1}^K\E_{\mathcal A} ( X_k X_{k+\ell})  \big )  - c^2 t \Big | 
 \leq  \sigma^2 \Big |  \frac{[nt]}{n} - t \Big |  + 
2  \E ( \zeta^2_{0})  \Big |  \frac{1}{n}\sum_{k=1}^{[nt]}   \sum_{\ell=1}^K{\bf 1}_{\phi_k \geq \ell } - t \sum_{i \geq 1} i \pi_i  \Big |  \\ + 2 \Big | \frac{1}{n}\sum_{k=1}^{[nt]} \sum_{\ell=1}^K \sum_{i=1}^{\ell}  \sum_{r \geq 0} U(k,\ell,i,r) -   \sum_{r \geq 1} \E ( \zeta_0 \zeta_r ) \sum_{j=1}^r (P^j)_{0,r-j}  \Big | =: I^{(1)}_n(t)+I^{(2)}_{n,K}(t)+I^{(3)}_{n,K}(t) \, ,
\end{multline}
where
\bbb \label{notationfortoprovealternativesecondcondp1-spet17}
U(k,\ell,i,r) := {\bf 1}_{\phi_k = \ell-i}{\bf 1}_{\phi_{k+\ell} = r} \E ( \zeta_{0} \zeta_{ r+i})\, .
\eee
Let us first prove that, for any $\varepsilon >0$  and any $t \in [0,1]$,
\bbb \label{toprovealternativesecondcondp2}
\lim_{m \rightarrow \infty} \limsup_{n \rightarrow \infty} \p \Big (  \frac{1}{m}  \sum_{K=1}^{m-1}  \big | I^{(2)}_{n,K}(t) \big |   \geq \varepsilon \Big )  =0  \, .
\eee
With this aim, note that since the Markov chain is irreducible and positive recurrent, by the law of large numbers for Markov chains, for any $t \in [0,1]$, 
\bbb \label{toprovealternativesecondcondp2(1)}
\lim_{n \rightarrow \infty}  \frac{1}{n}  \sum_{k=1}^{[nt]}  \sum_{\ell \geq 1}   {\bf 1}_{\phi_k \geq \ell } = \lim_{n \rightarrow \infty}  \frac{1}{n}  \sum_{k=1}^{[nt]} 
\phi_k =  t  \sum_{i \geq 1} i \pi_i \  \text{ a.s.} \, ,
\eee
whatever the initial law is.  In addition, by  Fact \ref{computationPrequalRWRTS}, 
\[
 \frac{1}{n}  \sum_{k=1}^{[nt]}  \sum_{\ell \geq  K+1}  \E_{\phi_0 =0}  ( {\bf 1}_{\phi_k \geq \ell } )  =  \frac{1}{n}  \sum_{k=1}^{[nt]}  \sum_{\ell \geq  K+1}  \p_{\phi_0 =0}  ( N_{k-1} = N_{k + \ell -1}) \leq  \sum_{\ell \geq  K+1}  \p ( \tau > \ell) \, .
\]
Hence 
\bbb  \label{toprovealternativesecondcondp2(2)}
 \frac{1}{m}  \sum_{K=1}^{m-1}  \frac{1}{n}  \sum_{k=1}^{[nt]}  \sum_{\ell \geq  K+1}  \E_{\phi_0 =0}  ( {\bf 1}_{\phi_k \geq \ell } )   \leq \frac{1}{m}  \sum_{K=1}^{m-1}  \sum_{\ell \geq  K+1}  \p ( \tau > \ell)  \rightarrow 0\, , \,  \text{as} \ m \rightarrow \infty \, . 
\eee
The convergence \eqref{toprovealternativesecondcondp2} follows from \eqref{toprovealternativesecondcondp2(1)} and \eqref{toprovealternativesecondcondp2(2)}.  To end the proof of 
\eqref{toprovealternativesecondcond} (and then of \eqref{condSqNSavecE0rho}), it remains to show that, for any $\varepsilon >0$  and any $t \in [0,1]$,
\bb 
\lim_{m \rightarrow \infty} \limsup_{n \rightarrow \infty} \p \Big (  \frac{1}{m}  \sum_{K=1}^{m-1}  \big | I^{(3)}_{n,K}(t) \big |   \geq \varepsilon \Big )  =0  \, .
\ee
This follows from tedious computations involving in particular again the ergodic Theorem for Markov chains and the fact that, as already mentioned, $ \sum_{m \geq 0} |\E ( \zeta_{0} \zeta_{ m}) | < \infty$ under the first part of condition \eqref{condARWRTSpart1stat}. The detailed proof is left to the reader.

\smallskip 

To end the proof of the corollary it remains to prove Facts \ref{computationPrequalRWRTS} and \ref{computationofthecondexpectationstat}. 

\smallskip

\noindent {\it Proof of Fact \ref{computationPrequalRWRTS}.} 
To prove the fact above we denote by $\tilde{\nu}$  the last renewal before time $k-m$, and write
\begin{align*}
\rP(N_{k-m} = N_k) 
& = \sum_{r=0}^{k-m}\rP(\tilde{\nu} = r, \tau_{\tilde{\nu} + 1} > k-r)
 = \sum_{r=0}^{k-m}\rP(\tilde{\nu} = r, \tau_{r + 1} > k-r)\\
& \leq \sum_{r=0}^{k-m}\rP(\tilde{\nu} = r, \tau_{r + 1} > m)
\leq \rP(\tau > m) \, .
\end{align*}
\qed

\medskip

\noindent {\it Proof of Fact \ref{computationofthecondexpectationstat}.} 
Note first that $Y_j = j + \phi_j= \nu_r$, for $\nu_{r-1}< j \leq \nu_r$, and by definition of $N_k$,
we have that $N_k - 1 = \max\{j \geq 1: \nu_j \leq k\}$, that is the last time the increasing random walk $\{\nu_j: j\geq 0\}$ is below the level $k$. 
Thus $N_k$ is the hitting time of the set $(k,\infty)$, or the time of first entry in the interval $[k+1,\infty)$ of the random walk $\{\nu_j: j\geq 0\}$. 
Thus $\nu_{N_{k}-1} < k+1$ and $\nu_{N_{k}} \geq k+1$. This fact with $k$ replaced by $k-1$ gives
$$\nu_{N_{k-1}-1} < k \leq \nu_{N_{k-1}}.$$
Since 
$Y_j = \nu_r$ for $\nu_{r-1}< j \leq \nu_r$ , we can conclude that $Y_k = \nu_{N_{k-1}}$ for all $k\geq 1$.
Then once again writing $\mathcal{A}$ for $\sigma(\phi_j, j\in \mathbb{Z})$, 
$\mathcal{F}_k = \sigma( \mathcal{A}, X_j, 1 \leq j\leq k)$, and $\rE_{\mathcal{A}}(\cdot)$ for $\rE(\cdot | \mathcal{A})$, we have
\[
\rE \big( X_k \big| \mathcal{F}_{k-m} \big)
= \rE_{\mathcal{A}} \big( \zeta_{Y_k} \big| \zeta_{Y_1},\dots, \zeta_{Y_{k-m}}\big)
= \rE_{\mathcal{A}} \big( \zeta_{\nu_{N_{k-1}}} \big| \zeta_{\nu_1},\dots, \zeta_{ \nu_{N_{k-m-1}}}\big) \, .
\]
Now, writing $k$ instead of $k-1$ to simplify notation, we have
\[
\left\| \rE_{\mathcal{A}}\big( \zeta_{\nu_{N_{k}}} \big| \zeta_{\nu_1},\dots, \zeta_{ \nu_{N_{k-m}}}  \big) \right\|_2 \leq
\left \| \rE_{\mathcal{A}}\big( \zeta_{\nu_{N_{k}}} \big| \mathcal{G}_{ \nu_{N_{k-m}}} \big) \right\|_2 
  \rE \big[ b^2(\nu_{N_{k}} - \nu_{N_{k-m}}) \big] \, .\]
Recall that $N_k$ is the first time the positive random walk $\nu_j$ exceeds the value $k$.
Thus 
\[\nu_{N_{k}} - \nu_{N_{k-m}} = \nu_{\sigma_k^{+}} - \nu_{\sigma_{k-m}^{+}}> k - (\nu_{\sigma_{k-m}^{+}} - (k-m)) - (k-m) = m - (\nu_{\sigma_{k-m}^{+}}- (k-m)) \, ,\]
where $\sigma_k^{+} = \inf\{j\geq 0 : \nu_j > k \}$, and thus $\nu_{\sigma_{k-m}^{+}}- (k-m)$ is the overshoot. A quick calculation gives for $l=[m/2]$:
\begin{align*}
\rE \big[ b^2(\nu_{N_{k}} - \nu_{N_{k-m}}) \big]
&= \rE\big[ b^2 (\nu_{N_{k}} - \nu_{N_{k-m}}) \mathbf{1}\{N_{k-m}\leq N_{k-l}-1\} \big] \\
&\quad + \rE\big[ b^2 (\nu_{N_{k}} - \nu_{N_{k-m}}) \mathbf{1}\{N_{k-m}=N_{k-l}\} \big]\\
&=: I_1 + I_2 \, ,
\end{align*}
where the first event corresponds to there being at least one renewal between times $k-m$ and $k-l$, while the second event corresponds to no renewals. Then, we have
\begin{align*}
I_1
&= \rE\big[ b^2 \left( \nu_{N_{k}} - \nu_{N_{k-l}-1 } + \nu_{N_{k-l}-1 } - \nu_{N_{k-m}} \right) \mathbf{1}\{N_{k-m}\leq N_{k-l}-1\} \big]\\
&\leq a^2  ([m/2]) \rP(N_{k-l} > N_{k-m}) \leq b^2([m/2]) \, .
\end{align*}
Indeed, by previous considerations, we have that 
$\nu_{N_{k}} - \nu_{N_{k-l}-1}> l$, while on the event $\{N_{k-m}\leq N_{k-l}-1\} $, clearly $\nu_{N_{k-l}-1} - \nu_{N_{k-m}} \geq 0$ and thus
$$\nu_{N_{k}} - \nu_{N_{k-l}-1 } + \nu_{N_{k-l}-1 } - \nu_{N_{k-m}} \geq l \, .$$
Therefore, since $b(\cdot)$ is non-increasing, 
$$b^2 \big(\nu_{N_{k}} - \nu_{N_{k-l}-1 } + \nu_{N_{k-l}-1 } - \nu_{N_{k-m}} \big) \leq b^2 (l) \, .$$
For the second term we have $I_2 \leq b^2(0) \rP(N_{k-m} = N_{k-l})$. So, according to Fact  \ref{computationPrequalRWRTS}, we get that for $k \geq m$,
\[
I_2 \leq b^2(0) \rP(\tau \geq  m - l ) \leq b^2 (0) \rP(\tau \geq [m/2]) \, .
\]
Taking into account the bounds of $I_1$ and $I_2$, the first inequality of Fact \ref{computationofthecondexpectationstat} follows. 

To prove the second one, we proceed similarly.  Since $Y_k = \zeta_{N_{k-1}} $,  for all $k \geq 0$, we have 
\begin{multline*}
\E ( X_{k+m} X_{k+m+\ell}| {\mathcal F}_k )  -\E_{\mathcal A} ( X_{k+m} X_{k+m+\ell})  \\
= \E_{\mathcal A} \big ( \zeta_{\nu_{N_{k+m-1}}} \zeta_{\nu_{N_{k+m + \ell-1}}}|  \zeta_{\nu_1}  , \ldots, \zeta_{\nu_{N_{k-1}} } \big )-\E_{\mathcal A} \big ( \zeta_{\nu_{N_{k+m-1}}} \zeta_{\nu_{N_{k+m + \ell-1}} }\big )  \, .
\end{multline*}
Therefore
\[
\Vert \E ( X_{k+m} X_{k+m+\ell}| {\mathcal F}_k )  -\E_{\mathcal A} ( X_{k+m} X_{k+m+\ell})  \Vert_{1 } \leq  \E  \big ( c ( \nu_{N_{k+m-1}} - \nu_{N_{k-1}} )  \big )   \, ,
\]
and the second inequality follows from the previous computations with $c(\cdot)$ replacing $b^2(\cdot)$. $\diamond$

\section{Appendix}

In this section we state a useful general lemma on subsequences.

\begin{lemma}
\label{lmaconvergence} Let $(A_{j}(n,m);n\geq1,m\geq1)$, $j=1,\ldots,J$, be a
double index family of nonnegative real numbers such that, for any
$j=1,\ldots,J$,
\[
\lim_{m\rightarrow\infty}\lim\sup_{n\rightarrow\infty}A_{j}(n,m)=0\,.
\]
Then for any $u_{n}\rightarrow\infty$, there exists a sequence of positive
integers $m_{n}$ such that $m_{n}\rightarrow\infty$, $m_{n}\leq u_{n}$ and,
for any $j=1,\ldots,J$,
\[
\lim_{n\rightarrow\infty}A_{j}(n,m_{n})=0.
\]

\end{lemma}

\textbf{Proof}. First, we observe that by considering the function
\[
A(n,m)=\max_{1\leq j\leq J}A_{j}(n,m) \, ,
\]
the lemma reduces to the case $J=1$.

Construct two strictly increasing positive integer sequences $\ell_{k}$ and
$n_{k}$ such that for all $n\geq n_{k}$,
\[
A(n,\ell_{k})\leq\frac{1}{k}\,.
\]
Let $g(n)=k$ for $n_{k}<n\leq n_{k+1}$ starting with $k=1$ and $g(n)=1$ for
$n\leq n_{1}$. Then, $g(n)$ is non-decreasing, $g(n)\rightarrow\infty$ and for
all $n$ such that $n_{k}<n\leq n_{k+1}$ we have
\[
n_{g(n)}=n_{k}<n \, .
\]
Now, let $G(n)$ be a sequence of positive integers such that $G(n)\leq g(n)$
and $G(n)\rightarrow\infty$ and $\ell_{G(n)}\leq u_{n}$. Then,
\[
n_{G(n)}\leq n_{g(n)}=n_{k}<n \, .
\]
Finally, let $m_{n}=\ell_{G(n)}$. Then, obviously $m_{n}\leq u_{n},$
$m_{n}\rightarrow\infty\;$and
\[
A(n,m_{n})=A(n,\ell_{G(n)})\leq\frac{1}{G(n)}\rightarrow0
\]
which proves the lemma. $\diamond$

\medskip

\noindent\textbf{Acknowledgements.} The research of Magda Peligrad was
partially supported in part by the NSF grant DMS-1512936 and by a Taft
Research Center award.


\begin{thebibliography}{99}                                                                                               %


\bibitem {Aldous78}Aldous, D. (1978). Stopping times and tightness.
\textit{Ann. Probab.} \textbf{6}, 335--340.

\bibitem {Bi1968}Billingsley, P. (1968). \textit{Convergence of Probability
Measures}, John Wiley \& Sons, New York.

\bibitem {Bi95}Billingsley, P. (1995). \textit{Probability and measure}. Third
edition. Wiley Series in Probability and Mathematical Statistics. A
Wiley-Interscience Publication. John Wiley \& Sons, Inc., New York.

\bibitem {Br}Bradley, R.C. (2007). \textit{Introduction to strong mixing
conditions}. 3 Volumes, Kendrick Press.

\bibitem {CDM}Cuny, C., Dedecker, J. and Merlev\`{e}de, F. (2017). Large and
moderate deviations for the left random walk on $GL_{d}({\mathbb{R}})$.
\textit{\ ALEA Lat. Am. J. Probab. Math. Stat.} {\textbf{1}}\textbf{4}, 503--527.

\bibitem {CM14}Cuny, C. and Merlev\`{e}de, F. (2014). On martingale
approximations and the quenched weak invariance principle. \textit{Ann.
Probab.} \textbf{\ 42}, 760--793.

\bibitem {De08}Dedecker, J. (2008). In\'egalit\'es de Hoeffding et
th\'eor\`eme limite central pour les fonctions peu r\'eguli\`eres de cha\^ines
de Markov non irr\'eductibles. \textit{\ num\'ero sp\'ecial des Annales de
l'ISUP} \textbf{\ 52}, 39--46.

\bibitem {DMP}Dedecker, J., Merlev\`{e}de, F and Peligrad, M. (2014). A
quenched weak invariance principle. \textit{Ann. Inst. Henri Poincar\'{e}
Probab. Stat.} \textbf{\ 50}, 872--898.

\bibitem {DMPU}Dedecker, J. , Merlev\`{e}de, F., Peligrad, M. and Utev, S.
(2009). Moderate deviations for stationary sequences of bounded random
variables. \textit{Ann. Inst. Henri Poincar\'{e} Probab. Stat.} \textbf{45}, 453--476.

\bibitem {Dob}Dobrushin, R. (1956). Central limit theorems for non-stationary
Markov chains I, II. \textit{Theory of Probab. and its Appl.} \textbf{1},
65--80, 329--383.

\bibitem {GanHa}G\"anssler, P. and H\"ausler, E. (1979). Remarks on the
functional central limit theorem for martingales. \textit{Z. Wahrsch. Verw.
Gebiete} \textbf{50}, no. 3, 237--243.

\bibitem {Gordin}Gordin, M. I. (1969). The central limit theorem for
stationary processes, \textit{Soviet. Math. Dokl.} \textbf{10}, 1174--1176.

\bibitem {GorPel}Gordin, M. and Peligrad, M. (2011). On the functional CLT via
martingale approximation. \textit{Bernoulli Journal} \textbf{17}, 424-440.

\bibitem {hh}Hall, P. and Heyde, C. C. (1980). \textit{Martingale limit theory
and its application.} Academic Press, New York-London.

\bibitem {He82}Helland, I. S. (1982). Central limit theorems for martingales
with discrete or continuous time. \textit{Scand. J. Statist.} \textbf{9}, no.
2, 79--94.

\bibitem {Heyde}Heyde, C. C. (1974). On the central limit theorem for
stationary processes. \textit{Z. Wahrsch. verw. Gebiete.} \textbf{30}, 315-320.

\bibitem {Krengel}Krengel, U. (1985). \textit{Ergodic theorems}, de Gruyter
Studies in Mathematics, 6. Walter de Gruyter \& Co., Berlin.

\bibitem {MW}Maxwell, M. and Woodroofe, M. (2000). Central limit theorem for
additive functionals of Markov chains. \textit{Ann. Probab.} \textbf{28}, 713--724.

\bibitem {McL75}McLeish, D. L. (1975). Invariance principles for dependent
variables. \textit{Z. Wahrscheinlichkeitstheorie und Verw. Gebiete}
\textbf{\ 32}, 165--178.

\bibitem {McL77}McLeish, D. L. (1977). On the invariance principle for
non-stationary mixingales. \textit{\ Ann. Probab.} \textbf{5}, 616--621.



\bibitem {PU1}Peligrad, M. and Utev, S. (2005). A new maximal inequality and
invariance principle for stationary sequences. \textit{Ann. Probab.}
\textbf{33}, 798--815.

\bibitem {PU2}Peligrad, M. and Utev, S. (2006). Central limit theorem for
stationary linear processes. \textit{Ann. Probab.} \textbf{34}, 1608--1622.

\bibitem {PU3}Peligrad, M. (2012). Central limit theorem for triangular arrays
of Non-Homogeneous Markov chains. \textit{Probability Theory and Related
Fields} \textbf{154}, 409--428.

\bibitem {PU3}Shao, Q. (1989). On the invariance principle for $\rho$-mixing
sequences of random variables. \textit{Chinese Ann. Math. (Ser. B)}
\textbf{10}, 427--433.

\bibitem {PU2}Sethuraman, S. and Varadhan, S. R. S. (2005). A martingale proof
of Dobrushin's theorem for non-homogeneous Markov chains. \textit{Electron. J.
Probab.} \textbf{10}, 1221--1235.

\bibitem {Utev90}Utev, S. (1990). Central limit theorem for dependent random
variables. \textit{Probab. Theory Math. Statist.} Vol. II (Vilnius, 1989),
519--528, ''Mokslas''', Vilnius.

\bibitem {Utev91}Utev, S. (1991). Sums of random variables with $\phi$-mixing
[translation of \textit{Trudy Inst. Mat.} (Novosibirsk) 13 (1989), Asimptot.
Analiz Raspred. Sluch. Protsess., 78-100]. \textit{Siberian Adv. Math.}
\textbf{1}, 24--155.

\bibitem {WuZh}Wu, W.B. and Zhao, Z. (2008). Moderate deviations for
stationary processes. \textit{Statistica Sinica} \textbf{18}, 769--782.

\bibitem {ZaWoo}Zhao, O. and Woodroofe, M. (2008). On Martingale
approximations, \textit{Annals of Applied Probability} \textbf{18}, 1831-1847.
\end{thebibliography}
\end{document}